\RequirePackage{fix-cm}
\documentclass[smallextended]{svjour3}       % onecolumn (second format)
\smartqed  % flush right qed marks, e.g. at end of proof

\usepackage{mathrsfs}
\usepackage{cite}
\usepackage[top=1in, bottom=1in, left=1.20in, right=1.20in]{geometry}
\usepackage{latexsym, amsmath, amssymb}
\usepackage{graphicx,epsfig}

\usepackage{amsfonts}
\usepackage{amscd}
\usepackage{amsbsy}
\usepackage{bm}
\usepackage{color}
\newcommand{\eqnref}[1]{(\ref {#1})}

\newcommand{\RR}{\mathbb{R}}

\newcommand{\beq}{\begin{equation}}
\newcommand{\eeq}{\end{equation}}

\begin{document}

\title{{\itshape Numerical methods for reconstruction of source term of heat equation from the final overdetermination}
\thanks{Fang is supported by NSF grants  No. NSFC70921001 and No. 71210003. Deng and Li are supported by NSF grants No. NSFC11301040,}
}


\author{Xiaoping Fang         \and
        Youjun Deng \and
        Jing Li %etc.
}


\institute{X. Fang \at
              Postdoctoral, Management Science and Engineering Postdoctoral Mobile Station, School of Business; School of Mathematics and Statistics, Central South University, Changsha, Hunan 410083, P. R. China. \\
              \email{fxpmath@csu.edu.cn}           %  \\
           \and
           Y. Deng \at
              Corresponding author. School of Mathematics and Statistics,
   Central South University, Changsha, Hunan 410083, P. R. China.
                \email{youjundeng@csu.edu.cn, dengyijun\_001@163.com}           %  \\
                \and
                J. Li \at
                Department of Mathematics, Changsha University of Science and Technology, Changsha, Hunan 410004, P.R.China.
                \email{lijingnew@126.com}
}

\date{Received: date / Accepted: date}

\maketitle

\begin{abstract}
This paper deals with the numerical methods for the reconstruction of source term in linear
parabolic equation from final overdetermination. We assume that the source term
has the form $f(x)h(t)$ and $h(t)$ is given, which guarantees the uniqueness of the inverse problem
of determining the source term $f(x)$ from final overdetermination. We present the regularization methods for 
reconstruction of the source term in the whole real line and with Neumann boundary
conditions. Moreover, we show the connection of the solutions between the problem with Neumann boundary conditions and 
the problem with no boundary condition (in the whole real line) by using extension method. 
Numerical experiments are done for the inverse problem with
the boundary conditions.\\
\noindent {\footnotesize {\bf Mathematics subject classification
(MSC2000):} 35R30, 35C20}\\
%{\em Keywords: Linear parabolic equation, source term, inverse problem, numerical methods}
\end{abstract}

\section{Introduction}
We consider the reconstruction of the source term in the following mathematical model
\begin{equation}
\label{eq:101}
\left\{
\begin{array}{ll}
 u_t = ku_{xx} + f(x)h(t) & 0<x<1, \ \ 0<t\leq T, \\
 u(x,0)=\mu_0(x) & x\in (0,l), \\
 u_x(0,t)= u_x(1,t)=0  & t \in (0,T],
\end{array}
\right.
\end{equation}
where $k$ is the heat conductivity, $f(x)$ and $h(t)$ relate to the source term. $\mu_0(x)$ is the initial status. $l$ and
$T$ are finite numbers.
If all those parameters are given then the direct problem (\ref{eq:101}) has
a unique solution. The inverse problem here is the determination of the source term $f(x)$
from the final state observation $\mu_T(x)=u(x,T)$.

The mathematical model (\ref{eq:101}) arises in various physical and engineering settings, in particular in
hydrology \cite{Bea:1972}, material sciences \cite{RHN:1987}, heat transfer \cite{Win:1997} and transport problems \cite{ZBe:1995}, etc.
The inverse problem in determination of source term has been studied intensively for decades
(cf., e.g., \cite{CDu:1980,CDu:1998,Iva:1998,LYH:2006,ODEK13}). The identification of an unknown state-dependent source term in a reaction-diffusion
equation is considered in  \cite{CDu:1980,CDu:1998}. In \cite{Iva:1998} the uniqueness of the inverse source problem with
arbitrary boundary conditions has been proved under several additional conditions. In \cite{LYH:2006}, the source function $f(x)$ is assumed
to be the sum of a known function $f(x)=\sum_{i=1}^I \rho(x-a_i)$ with $I$ different locations and the locations $a_i$ are determined by three non-collinear measurement points.
On the other hand, the inverse problems for parabolic equations with
final overdetermination also have been considered by lots of authors (see \cite{IBu:1995,Cho:1994,CYm:1996,Has:2007} and the
references there in). However, numerical methods for uniquely solving the inverse source term $f(x)$ in (\ref{eq:101}) without using
further data concerning the solution $u(x,t)$ are seldom. We shall provide the numerical solution for solving the inverse source
problem (\ref{eq:101}) and more importantly, we show the relationship between solution of the boundary problem (\ref{eq:101}) and its
corresponding no boundary problem (in the whole real line). 
We only consider the one dimensional problem here to simplify our calculation and point out the
main idea. The method can in fact be implemented in two dimensional or higher dimensional problem.

In this paper, we shall first consider the heat conduction problem in the whole real line, which is
\begin{equation}
\label{eq:103}
\left\{
\begin{array}{ll}
 u_t = ku_{xx} + f(x)h(t) & x\in \RR, \ \ 0<t\leq T, \\
 u(x,0)=\mu_0(x) & x\in \RR. \\
\end{array}
\right.
\end{equation}
where we suppose $f(x),\mu_0(x)\in L^2(\RR)$ and $h(t)\in L^2(0,T)$.
It is easy to see that the solution of (\ref{eq:103}) (cf. \cite{Eva:1998}) is
\begin{equation}
\label{eq:001}
u(x,t)  =  \int_{-\infty}^{\infty} \frac{1}{\sqrt{4k\pi t}}e^{-\frac{(x-y)^2}{4kt}} \mu_0(y) dy
+ \int_{0}^{t}\int_{-\infty}^{\infty} \frac{1}{\sqrt{4k\pi (t-s)}}e^{-\frac{(x-y)^2}{4k(t-s)}} f(y)h(s) dyds.
\end{equation}
By taking the Fourier transform with respect
to $x$ we can immediately get 
$$\hat{u}_t(\xi,t)= -k\xi^2 \hat{u}(\xi,t) + \hat{f}(\xi)h(t)$$ 
and by initial 
condition in \eqnref{eq:103} there holds
\begin{equation}
\label{eq:104}
\hat{u}(\xi,t)=\hat{\mu}_0(\xi)e^{-k\xi^2t}+\hat{f}(\xi)\int_0^th(s)e^{-k\xi^2(t-s)}ds .
\end{equation}
For the inverse problem with $\mu_T(x)$ being measured, from (\ref{eq:104}) we have
$$
\hat{\mu}_T(\xi)=\hat{\mu}_0(\xi)e^{-k\xi^2T}+\hat{f}(\xi)\int_0^Th(s)e^{-k\xi^2(T-s)}ds,
$$
and so
\begin{equation}
\label{eq:106}
\hat{f}(\xi)=\frac{\hat{\mu}_T(\xi)-\hat{\mu}_0(\xi)e^{-k\xi^2T}}{\int_0^Th(s)e^{-k\xi^2(T-s)}ds}.
\end{equation}
Furthermore, the solution $u(x,t)$ can also be written as
\begin{eqnarray}
\hat{u}(\xi,t) & = & \hat{\mu}_0(\xi)e^{-k\xi^2t}+(\hat{\mu}_T(\xi)-\hat{\mu}_0(\xi)e^{-k\xi^2T})
\frac{\int_0^th(s)e^{-k\xi^2(t-s)}ds}{\int_0^Th(s)e^{-k\xi^2(T-s)}ds}  \nonumber \\
& = & \hat{\mu}_0(\xi)e^{-k\xi^2t}\frac{\int_t^{T}h(s)e^{k\xi^2s}ds}{\int_0^Th(s)e^{k\xi^2s}ds}
+ \hat{\mu}_T(\xi)\frac{\int_0^th(s)e^{-k\xi^2(t-s)}ds}{\int_0^Th(s)e^{-k\xi^2(T-s)}ds}.
\label{eq:107}
\end{eqnarray}
The relation (\ref{eq:106}) tells that if $h(t)$ is appropriately given in $[0,T]$ such that the
 denominator in (\ref{eq:106}) is nonzero for every $\xi$ (or be nonzero in the distribution meaning),
  and $\mu_0(x)\in L^2(\RR)$ and $\mu_T(x)\in L^2(\RR)$ are given for $x\in \RR$, then
$\hat{f}(\xi)$ and so $f(x)$ can be reconstructed uniquely. The argument is also suitable for (\ref{eq:101}) as we
shall see.

In this paper, we consider the determination of source term $f(x)$ in both \eqnref{eq:101} and
\eqnref{eq:103}. For the sake of simplicity, we assume that
$h(t)$ is identically non-positive or non-negative function in $[0,T]$ and $C_h:=\int_0^T|h(s)|ds>0$.
This assumption is quite nature. For
example, in the case where the heat is provided by a single kind of radioactive isotope, we can set
$h(t)=e^{-\lambda t}$ with a constant $\lambda>0$ (cf. \cite{LYH:2006}). Furthermore, we suppose
$f(x)\in H^p(\RR)$, $p\geq 0$, where $\|f\|_p$ is defined as the norm of $f(x)$ in $H^p(\RR)$
\begin{equation}
\label{eq:108}
\|f\|_p:=(\int_{-\infty}^{\infty}(1+{\xi}^2)^p|\hat{f}(\xi)|^2d\xi)^{1/2}.
\end{equation}
This paper is organized as follows. In section 2, we analyze the severe ill-posedness of the
reconstruction of the source term in \eqnref{eq:103}, then
we introduce the iterative method to solve the inverse problem. Convergence rates are given
under both {\em a priori} and  {\em a posteriori} stopping rules. In section 3, we use 
the same fundamental solution method to get the solution of \eqnref{eq:101}, by extending the source term and the initial state to
the whole real axis and then show the solution of the Neumann boundary problem is actually 
another form of the solution to \eqnref{eq:101} by separating variables method. In fact,  the solution to any kind of boundary problem
can be got by the method extending the source and initial terms to the whole region other than by separating variables method. We
only use the Neumann boundary problem as an example.
We then give the frequency cut-off technique to solve the inverse problem. Numerical experiments
for boundary problem are done in section 4 and show attractive results.

\section{No boundary restriction case}
In this section, we discuss the reconstruction of $f(x)$ in $L^2(\RR)$. We shall keep the denotation of norm in $L^2(\RR)$
as $\|\cdot\|$. The direct problem is
\begin{equation}
\label{eq:201}
\left\{
\begin{array}{ll}
 u_t = ku_{xx} + f(x)h(t) & (x,t) \in \RR\times (0,T], \\
 u(x,0)=\mu_0(x) & x\in \RR, \\
\end{array}
\right.
\end{equation}
and the inverse problem is reconstruction of $f(x)$ with final overdetermination $u(x,T)=\mu_T(x)$.
We have already got $\hat{f}(\xi)$ in (\ref{eq:106}). However, it is a severely ill-posed problem.
In fact, the denominator in (\ref{eq:106}) decrease to zero exponentially as $\xi \rightarrow \infty$.
Thus small perturbations in the measured data $\mu_T(x)$ may produce high frequency parts in $\hat{f}(\xi)$
and make the reconstruction quite unstable. %On the other hand, we can first solve $u(x,t)$ by (\ref{eq:107}) and
%then using (\ref{eq:104}) to solve $f(x)$. Since $f(x)$ can be stably reconstructed by using (\ref{eq:104}),
%We only need to think about the reconstruction of $u(x,t)$. From (\ref{eq:107}) we know that the first term
%in the right hand of the equation is stable, while the second term is severely ill-posed since we can
%rewrite it as
%$$\hat{u}(\xi,t)= \hat{\mu}_0(\xi)e^{-k\xi^2t}\frac{\int_t^{T}h(s)e^{k\xi^2s}ds}{\int_0^Th(s)e^{k\xi^2s}ds}
%+ \hat{\mu}_T(\xi)e^{k\xi^2(T-t)}\frac{\int_0^th(s)e^{k\xi^2s}ds}{\int_0^Th(s)e^{k\xi^2s}ds}$$
%and we see that the fraction part is less than one for all $\xi \in \RR$.
%In this section, we mainly discuss the first method, that is, first solve $f(x)$ by (\ref{eq:106}) and then
%solve $u(x,t)$ by (\ref{eq:104}).
Suppose that the measured
final overspecified data $\mu_T^{\delta}(x)$ satisfies
\begin{equation}
\label{eq:202}
\|\mu_T^{\delta}-\mu_T\|\leq \delta.
\end{equation}
We do not assume measurement error in $\mu_0(x)$ because the
perturbation of $\mu_0(x)$ will affect very little in the reconstruction of $f(x)$ and $u(x,t)$.
Denote $v(\xi,t)=e^{-k\xi^2t}$ and rewrite
(\ref{eq:106}) as
\begin{equation}
\label{eq:203}
\hat{g}(\xi):=\hat{\mu}_0(\xi)+\hat{f}(\xi)\int_0^Th(s)e^{k\xi^2s}ds=\hat{\mu}_T(\xi)e^{k\xi^2T}
=\frac{\hat{\mu}_T(\xi)}{v(\xi,T)}.
\end{equation}
We can now introduce similar iterative method to \cite{DYJ:2009,DYJ:2010} for solving $\hat{g}(\xi)$
\begin{equation}
\label{eq:204}
\hat{g}_n^{\delta}(\xi)=(1-\lambda)\hat{g}_{n-1}^{\delta}(\xi)+
\frac{\lambda}{v(\xi,T)}\chi_{\vartheta_1}\hat{\mu}_T^{\delta}(\xi)+\lambda(1-\chi_{\vartheta_1})\hat{\mu}_0(\xi),
\end{equation}
where $\lambda=\sqrt[N]{v(\xi,T)}<1$, $N$ is a nature number, $\chi_{\vartheta_1}$ denotes the characteristic function of interval
$[-\vartheta_1,\vartheta_1]$ and $\vartheta_1$ is large enough.
%%%First we give the following elementary lemma.
%%%\begin{lemma}
%%%Let $C_h=|\int_0^T h(s) ds|>0$ and $m(\xi)=|\int_0^T h(s) e^{k\xi^2 s} ds|$, where $h(t)$ is assumed to be
%%%identically non-positive or non-negative in $[0,T]$, then
%%%$$m(\xi) \geq C_h,$$
%%%and
%%%$$m(\xi)\geq C \frac{e^{k\xi^2\sigma}-1}{k\xi^2},$$
%%%where $C$ and $\sigma$ are dependent on $h(t)$.
%%%\end{lemma}
%%%\textbf{Proof.} It is sufficient to prove the second inequality. Since $C_h>0$, by the definition of the integral, there
%%%exists a segment, say $[t_1,t_2]$, and $t_2-t_1\geq \sigma$ such that $|h(t)|>\epsilon$ in it.
%%%\begin{eqnarray*}
%%%m(\xi)& = & \left|\int_0^T h(s) e^{k\xi^2 s} ds\right|=\int_0^T |h(s)| e^{k\xi^2 s} ds \\
%%%& \geq & \epsilon \int_{t1}^{t2} e^{k\xi^2 s} ds = \epsilon e^{k\xi^2t_1}\frac{e^{k\xi^2(t_2-t_1)}-1}{k\xi^2} \\
%%%& \geq & \epsilon \frac{e^{k\xi^2\sigma}-1}{k\xi^2}
%%%\end{eqnarray*}
%%%which completes the proof. \hspace{\fill} $\Box$
\subsection{A priori stopping rule}
By using the a priori stopping rule we have the following convergence theorem
\begin{theorem}
\label{Th:201} Let $u(x,t)$ be the exact temperature history of (\ref{eq:201}), $h(t)\not\equiv 0$ is identically
nonpositive or nonnegative in $[0,T]$ and $\mu_T^{\delta}(x)$ be the measured final temperature
satisfying (\ref{eq:202}). $f(x)$ satisfies $\|f\|_p\leq M$. Let $\hat{g}_k^{\delta}(\xi)$ be the
$k$-th iteration solution defined by (\ref{eq:204}) with $\hat{g}_0^{\delta}(\xi)=\hat{\mu}_0(\xi)$,
where $\vartheta_1\sim\sqrt{\frac{1}{(1+\sigma)kT}\left[ln\frac{M}{\delta}(ln\frac{M}{\delta})^{-\frac{1+\sigma}{2}p}\right]}$
, $\sigma\geq 0$. Suppose $u_n^{\delta}(x,t)$, $f_n^{\delta}(x)$
are solved via (\ref{eq:104}), (\ref{eq:106}) and the inverse Fourier transform for every $\hat{g}_k^{\delta}$,
respectively. If we choose $n\sim\lfloor\sqrt[N]{\frac{M}{\delta}}\rfloor$, then there holds
\begin{equation}
\label{eq:th201} {\| f_{n}^{\delta }-f
\|}^{2} \leq C(ln\frac{M}{\delta})^{-p}
\left(  M^{\frac{2}{1+\sigma}}\delta^{\frac{2\sigma}{1+\sigma}}+ M^{2}(\frac{ln\frac{M}{\delta}}
{ln\frac{M}{\delta}(ln\frac{M}{\delta})^{-\frac{1+\sigma}{2}p}})^{p} \right),
\end{equation}
and
\begin{equation}
\label{eq:th202}
{\| u_{n}^{\delta }(\cdot,t)-u(\cdot,t)\|}^{2}
\leq C_h(t)(ln\frac{M}{\delta})^{-p}
\left(  M^{\frac{2}{1+\sigma}}\delta^{\frac{2\sigma}{1+\sigma}}+ M^{2}(\frac{ln\frac{M}{\delta}}
{ln\frac{M}{\delta}(ln\frac{M}{\delta})^{-\frac{1+\sigma}{2}p}})^{p} \right),
\end{equation}
for $\delta \rightarrow 0$,
where $C$ is a constant independent of $\delta$ and $M$ and $C_h(t)=C(\int_0^t |h(s)|ds)^2$.
\end{theorem}
\textbf{Proof.} By the iteration (\ref{eq:204}) we have
\begin{eqnarray*}
\hat{g}_n^{\delta}(\xi) & = & (1-\lambda)\hat{g}_{n-1}^{\delta}(\xi)+
\frac{\lambda}{v(\xi,T)}\chi_{\vartheta_1}\hat{\mu}_T^{\delta}(\xi)+ \lambda(1-\chi_{\vartheta_1})\hat{\mu}_0(\xi)\\
& = & (1-\lambda)^n\hat{\mu}_0(\xi)+\sum_{i=0}^{n-1} (1-\lambda)^i\left[\frac{\lambda}{v(\xi,T)}
\chi_{\vartheta_1}\hat{\mu}_T^{\delta}(\xi)+ \lambda(1-\chi_{\vartheta_1})\hat{\mu}_0(\xi)\right] \\
& = & (1-\lambda)^n\hat{\mu}_0(\xi) + \sum_{i=0}^{n-1} (1-\lambda)^i\lambda(1-\chi_{\vartheta_1})\hat{\mu}_0(\xi)
+\sum_{i=0}^{n-1} (1-\lambda)^i\frac{\lambda}{v(\xi,T)}\chi_{\vartheta_1}\hat{\mu}_T^{\delta}(\xi).
\end{eqnarray*}
Set $p_n(\lambda)=\sum_{i=0}^{n-1} (1-\lambda)^i$, $r_n(\lambda)=1-\lambda p_n(\lambda)=(1-\lambda)^n$, we have the elementary
results (cf. \cite{VG:1986})
\begin{eqnarray*}
 & & p_{n}(\lambda)\lambda^{\mu}\leq n^{1-\mu},
\mbox{for} \ \ \mbox{all} \ \ 0\leq \mu\leq 1, \\
& & r_{n}(\lambda)\lambda^{v}\leq \theta_v
(n+1)^{-v},
\end{eqnarray*}
where
$$\theta_v= \left \{
\begin{array}{ll}
1, & 0 \leq v \leq 1 \\
v^v, & v > 1.
\end{array}
\right. $$
\begin{eqnarray*}
\hat{g}_n^{\delta}(\xi)-\hat{g}(\xi)& = & r_n(\lambda)\hat{\mu}_0(\xi)+p_n(\lambda)\lambda (1-\chi_{\vartheta_1})
\hat{\mu}_0(\xi)+\frac{p_n(\lambda)\lambda}{v(\xi,T)}\chi_{\vartheta_1}\hat{\mu}_T^{\delta}(\xi)-\hat{g}(\xi) \\
& = & \frac{p_n(\lambda)\lambda}{v(\xi,T)}[\chi_{\vartheta_1}\hat{\mu}_T^{\delta}(\xi)-\hat{\mu}_T(\xi)
+v(\xi,T)(1-\chi_{\vartheta_1})\hat{\mu}_0(\xi)] - r_n(\lambda)[\hat{g}(\xi)-\hat{\mu}_0(\xi)] \\
& = & \frac{p_n(\lambda)\lambda}{v(\xi,T)}[\chi_{\vartheta_1}(\hat{\mu}_T^{\delta}(\xi)-\hat{\mu}_T(\xi))
+ v(\xi,T)(1-\chi_{\vartheta_1})\hat{f}(\xi)\int_0^Th(s)e^{k\xi^2s}ds] \\
& & - r_n(\lambda)\hat{f}(\xi)\int_0^Th(s)e^{k\xi^2s}ds.
\end{eqnarray*}
Thus
\begin{eqnarray*}
{\| \hat{f}_{n}^{\delta }-\hat{f}
\|}^{2} & = & \int_{-\infty}^{\infty} \left(\frac{ \hat{g}^{\delta}_n(\xi)-\hat{g}(\xi)}{\int_0^Th(s)e^{k\xi^2s}ds}\right)^2d\xi \\
& \leq & 2\int_{-\infty}^{\infty}\left(\frac{ \frac{p_n(\lambda)\lambda}{v(\xi,T)}[\chi_{\vartheta_1}(\hat{\mu}_T^{\delta}(\xi)-\hat{\mu}_T(\xi))
+ v(\xi,T)(1-\chi_{\vartheta_1})\hat{f}(\xi)\int_0^Th(s)e^{k\xi^2s}ds]}{\int_0^Th(s)e^{k\xi^2s}ds}\right)^2 d\xi\\
& & +2\int_{-\infty}^{\infty}\left(r_n(\lambda)\hat{f}(\xi)\right)^2d\xi :=  2 I_1 +2 I_2.
\end{eqnarray*}
Next, we give separated evaluation for $I_1$ and $I_2$. We have
\begin{eqnarray*}
I_1 & = & \int_{-\infty}^{\infty}\left(\frac{ \frac{p_n(\lambda)\lambda}{v(\xi,T)}[\chi_{\vartheta_1}(\hat{\mu}_T^{\delta}(\xi)-\hat{\mu}_T(\xi))
+ v(\xi,T)(1-\chi_{\vartheta_1})\hat{f}(\xi)\int_0^Th(s)e^{k\xi^2s}ds]}{\int_0^Th(s)e^{k\xi^2s}ds}\right)^2d\xi \\
& = & \int_{|\xi|\leq\vartheta_1} \frac{\left( \frac{p_n(\lambda)\lambda}{v(\xi,T)}\left[\chi_{\vartheta_1}
(\hat{\mu}^{\delta}_T(\xi)-\hat{\mu}_T(\xi))\right]\right)^2}{(\int_0^Th(s)e^{k\xi^2s}ds)^2}d\xi
 +\int_{|\xi|>\vartheta_1} (p_n(\lambda)\lambda\hat{f}(\xi))^2d\xi \\
& \leq & \frac{1}{C_h}e^{2k\vartheta_1^2T}\delta^2 +  \int_{|\xi|>\vartheta_1}\hat{f}(\xi)^2 d\xi \\
& \leq & \frac{1}{C_h^2} M^{\frac{2}{1+\sigma}} \delta^{\frac{2\sigma}{1+\sigma}}(ln\frac{M}{\delta})^{-p}
+\vartheta_1^{-2p}M^2 \\
& \leq & \frac{1}{C_h^2} M^{\frac{2}{1+\sigma}} \delta^{\frac{2\sigma}{1+\sigma}}(ln\frac{M}{\delta})^{-p}
+ CM^2\left[ln\frac{M}{\delta}(ln\frac{M}{\delta})^{-\frac{1+\sigma}{2}p}\right]^{-p} \\
& \leq & (ln\frac{M}{\delta})^{-p}
\left( \frac{1}{C_h^2} M^{\frac{2}{1+\sigma}} \delta^{\frac{2\sigma}{1+\sigma}}+ CM^{2}(\frac{ln\frac{M}{\delta}}
{ln\frac{M}{\delta}(ln\frac{M}{\delta})^{-\frac{1+\sigma}{2}p}})^{p} \right)
\end{eqnarray*}
where $C$ is the general constant depending on $k$ and $T$.
\begin{eqnarray*}
I_2 & = &\int_{-\infty}^{\infty}(r_n(\lambda)\hat{f}(\xi))^2d\xi \\
& = & \int_{|\xi|\leq\vartheta_1}(r_n(\lambda)\hat{f}(\xi))^2d\xi
+ \int_{|\xi|>\vartheta_1}(r_n(\lambda)\hat{f}(\xi))^2d\xi \\
& \leq & \int_{|\xi|\leq\vartheta_1} (r_n(\lambda)\lambda^N \frac{\hat{f}(\xi)}{v(\xi,T)})^2 d\xi
+ \vartheta_1^{-2p}M^2 \\
& \leq & N^{2N} (n+1)^{-2N}e^{2k\vartheta_1^2T}M^2+\vartheta_1^{-2p}M^2 \\
& \leq & C N^{2N} M^{\frac{2}{1+\sigma}} \delta^{\frac{2\sigma}{1+\sigma}} (ln\frac{M}{\delta})^{-p}
+ CM^2\left[ln\frac{M}{\delta}(ln\frac{M}{\delta})^{-\frac{1+\sigma}{2}p}\right]^{-p} \\
& \leq & (ln\frac{M}{\delta})^{-p}
\left( C N^{2N} M^{\frac{2}{1+\sigma}} \delta^{\frac{2\sigma}{1+\sigma}}+ CM^{2}(\frac{ln\frac{M}{\delta}}
{ln\frac{M}{\delta}(ln\frac{M}{\delta})^{-\frac{1+\sigma}{2}p}})^{p} \right).
\end{eqnarray*}
Finally we have
$${\| \hat{f}_{n}^{\delta }-\hat{f}
\|}^{2} \leq C(ln\frac{M}{\delta})^{-p}
\left(  M^{\frac{2}{1+\sigma}} \delta^{\frac{2\sigma}{1+\sigma}}+ M^{2}(\frac{ln\frac{M}{\delta}}
{ln\frac{M}{\delta}(ln\frac{M}{\delta})^{-\frac{1+\sigma}{2}p}})^{p} \right)$$
and
\begin{eqnarray*}
{\| \hat{u}_{n}^{\delta }(\cdot,t)-\hat{u}(\cdot,t)
\|}^{2} & = & \int_{-\infty}^{\infty}\left[(\hat{f}_{n}^{\delta }(\xi)-\hat{f}(\xi))\int_0^th(s)e^{-k\xi^2(t-s)}ds\right]^{2}d\xi \\
& \leq & (\int_0^t|h(s)|ds)^2(ln\frac{M}{\delta})^{-p}
\left(  M^{\frac{2}{1+\sigma}} \delta^{\frac{2\sigma}{1+\sigma}}+ M^{2}(\frac{ln\frac{M}{\delta}}
{ln\frac{M}{\delta}(ln\frac{M}{\delta})^{-\frac{1+\sigma}{2}p}})^{p} \right).
\end{eqnarray*}
The proof is completed by using the Parseval equality. \hspace{\fill} $\Box$
\begin{remark}
We see in Theorem 1 that if $p=0$ then the second term in (\ref{eq:th201}) and (\ref{eq:th202}) is just a
bounded term and does not converge when $\delta \rightarrow 0$. However, the fact that the second term turns to
zero when $\delta \rightarrow 0$ if due to that $\int_{|\xi|>\vartheta_1}(r_n(\lambda)\hat{f}(\xi))^2d\xi$ definitely turns
to zero (since $\vartheta_1 \rightarrow \infty$ and $f \in L^2(\RR)$). Thus if $p=0$, which means $f(x) \in L^2(\RR)$ and
$f(x) \not\in H^p(\RR), p>0$, and by choosing $\sigma>0$ then we only obtain the convergence of the solution but with no convergence rate.
\end{remark}

\subsection{A posteriori stopping rule}
We introduce the widely-used "discrepancy principle" due to Morozov
\cite{MO:1966} in the following
form:
\begin{equation}
\label{eq:221}
\|\frac{\hat{\mu}_T^{\delta}-v(\cdot,T)\hat{g}_{n_*}^{\delta}}{\varphi}\|
\leq \tau \delta^{\frac{1}{1+\sigma}} <
\|\frac{\hat{\mu}_T^{\delta}-v(\cdot,T)\hat{g}_{n_*}^{\delta}}{\varphi}\|,
\ \ \mbox{for} \ \ 0\leq n < n_*,
\end{equation}
where $\varphi(\xi)=\int_0^T h(s) e^{k\xi s}ds$ and $n_*$ is the first iteration step which satisfies the left
inequality of (\ref{eq:221}). With the discrepancy principle, we have similar convergence results
\begin{theorem}
\label{Th:201} Let $u(x,t)$ be the exact temperature history of (\ref{eq:201}), $h(t)\not\equiv 0$ is identically
nonpositive or nonnegative in $[0,T]$ and $\mu_T^{\delta}(x)$ be the measured final temperature
satisfying (\ref{eq:202}). $f(x)$ satisfies $\|f\|_p\leq M$. Let $\hat{g}_k^{\delta}(\xi)$ denote the
iterates defined by (\ref{eq:204}) with $\hat{g}_0^{\delta}(\xi)=\hat{\mu}_0(\xi)$,
where $\vartheta_1\sim\sqrt{\frac{1}{(1+\sigma)kT}\left[ln\frac{M}{\delta}(ln\frac{M}{\delta})^{-\frac{1+\sigma}{2}p}\right]}$
. Suppose $u_n^{\delta}(x,t)$, $f_n^{\delta}(x)$
are solved via (\ref{eq:104}), (\ref{eq:106}) and the inverse Fourier transform for every $\hat{g}_k^{\delta}$,
respectively. If we select (\ref{eq:221}) as the a posteriori stopping rule, then there holds
\begin{equation}
\label{eq:th221} {\| f_{n}^{\delta }-f
\|}^{2} \leq C(ln\frac{M}{\delta})^{-p}
\left(  M^{\frac{2}{1+\sigma}} \delta^{\frac{2\sigma}{1+\sigma}}+ M^{2}(\frac{ln\frac{M}{\delta}}
{ln\frac{M}{\delta}(ln\frac{M}{\delta})^{-\frac{1+\sigma}{2}p}})^{p} \right),
\end{equation}
and
\begin{equation}
\label{eq:th222}
{\| \hat{u}_{n}^{\delta }(\cdot,t)-\hat{u}(\cdot,t)\|}^{2}
\leq C_h(t)(ln\frac{M}{\delta})^{-p}
\left(  M^{\frac{2}{1+\sigma}} \delta^{\frac{2\sigma}{1+\sigma}}+ M^{2}(\frac{ln\frac{M}{\delta}}
{ln\frac{M}{\delta}(ln\frac{M}{\delta})^{-\frac{1+\sigma}{2}p}})^{p} \right),
\end{equation}
for $\delta \rightarrow 0$,
where $C$ is a constant independent of $\delta$ and $M$ and $C_h(t)=C(\int_0^t |h(s)|ds)^2$.
\end{theorem}
\textbf{Proof.}
It suffices to prove that $n_* \sim \lfloor \sqrt[N]{\frac{M}{\delta}} \rfloor$. In fact, by (\ref{eq:204}) we have
\begin{eqnarray*}
\|\frac{v(\cdot,T)\hat{g}_{n}^{\delta}-\hat{\mu}_T^{\delta}}{\varphi}\| & = &
\|\frac{(1-\lambda)[v(\cdot,T)\hat{g}_{n-1}^{\delta}-\chi_{\vartheta_1}\hat{\mu}_T^{\delta}]+\lambda (1-\chi_{\vartheta_1})v(\cdot,T)
\hat{\mu}_0+(1-\chi_{\vartheta_1})\hat{\mu}_T^{\delta}}{\varphi}\| \\
& = & \|\frac{r_n(\lambda)[v(\cdot,T)\hat{\mu}_{0}-\chi_{\vartheta_1}\hat{\mu}_T^{\delta}]+
p_n(\lambda)\lambda (1-\chi_{\vartheta_1})v(\cdot,T)\hat{\mu}_0+ (1-\chi_{\vartheta_1})\hat{\mu}_T^{\delta}}{\varphi}\| \\
& \leq & \|\frac{r_n(\lambda)\chi_{\vartheta_1}(\hat{\mu}_{T}-\hat{\mu}_T^{\delta})
+ (1-\chi_{\vartheta_1})(\hat{\mu}_T-\hat{\mu}_T^{\delta})}{\varphi}\|\\
& & +\|\frac{r_n(\lambda)\chi_{\vartheta_1}(\hat{\mu}_{T}-v(\cdot,T)\hat{\mu}_{0})+
(1-\chi_{\vartheta_1})(\hat{\mu}_T-v(\cdot,T)\hat{\mu}_0)}{\varphi}\|\\
& \leq & \frac{\delta}{C_h}  + \|r_n(\lambda)\chi_{\vartheta_1}v(\cdot,T)\hat{f}+
(1-\chi_{\vartheta_1})v(\cdot,T)\hat{f}\|\\
& \leq & \frac{\delta}{C_h} + N^N (n+1)^{-N} M +e^{-k\vartheta_1^2T}\vartheta_1^{-p}M \\
& \leq & \frac{\delta}{C_h} + N^N (n+1)^{-N} M + CM^{\frac{\sigma}{1+\sigma}} \delta^{\frac{1}{1+\sigma}}
(\frac{ln\frac{M}{\delta}}{ln\frac{M}{\delta}(ln\frac{M}{\delta})^{-\frac{1+\sigma}{2}p}})^{p} \\
& = &  \delta^{\frac{1}{1+\sigma}}(C_1M^{\frac{\sigma}{1+\sigma}}+\frac{1}{C_h}\delta^{\frac{\sigma}{1+\sigma}})
+ N^N (n+1)^{-N} M,
\end{eqnarray*}
where $C_1=(\frac{ln\frac{M}{\delta}}{ln\frac{M}{\delta}(ln\frac{M}{\delta})^{-\frac{1+\sigma}{2}p}})^{p}$ is a
bounded term. We come to the conclusion if we choose $\tau=C_1M^{\frac{\sigma}{1+\sigma}}+\frac{C_h+1}{C_h}\delta^{\frac{\sigma}{1+\sigma}}
$.

\section{Boundary condition case}
In this section, we consider the reconstruction of source term $f(x)$ in (\ref{eq:101}) with final
overdetermination $\mu_T(x)$. For the sake of simplicity we still use the sign $\|\cdot\|$ for norm in $L^2[0,1]$,
which we expect it would not cause confusion to the reader.
At the beginning, we analyze the solution of direct problem (\ref{eq:101}). We can actually use
the extension method to get the solution of (\ref{eq:101}). Specifically, it has the solution as follows
\begin{equation}
\label{eq:301}
u(x,t)= \int_{-\infty}^{\infty} \frac{1}{\sqrt{4k\pi t}}e^{-\frac{(x-y)^2}{4kt}} \mu_0(y) dy
+ \int_{0}^{t}\int_{-\infty}^{\infty} \frac{1}{\sqrt{4k\pi (t-s)}}e^{-\frac{(x-y)^2}{4k(t-s)}} f(y) h(s)dyds,
\end{equation}
where $\mu_0(x)$ and $f(x)$ are extended to the following form
\begin{eqnarray*}
\mu_0(x) & = & \int_0^1 \mu_0(y)  dy + 2\sum_{m=1}^{\infty} \cos(m\pi x) \int_0^1 \mu_0(y) \cos(m\pi y) dy, \\
f(x) & = & \int_0^1 f(y) dy +2\sum_{m=1}^{\infty} \cos(m\pi x) \int_0^1 f(y) \cos(m\pi y) dy.
\end{eqnarray*}
It is easy to verify that under the above extensions of $\mu_0(x)$ and $f(x)$, the boundary conditions are satisfied.
Furthermore the solution (\ref{eq:301}) is similar to (\ref{eq:001}). But we can not directly use the
iterative method mentioned in section 2 to solve the inverse problem to get $f(x)$ and $u(x,t)$. In order to use the
iterative method (\ref{eq:204}) to solve the inverse problem we need further extend $\mu_T(x)$ as
$$\mu_T(x)  =  \int_0^1 \mu_T(y)dy + 2\sum_{m=1}^{\infty} \cos(m\pi x) \int_0^1 \mu_T(y) \cos(m\pi y) dy.$$
However, the extended functions $\mu_0(x)$ and $\mu_T(x)$ are no longer $L^2$ integrable functions in $\RR$.
Thus some further analysis of (\ref{eq:301}) is needed.

In fact (\ref{eq:301}) can be changed into another solution form. To explain this, firstly we give the following lemma.
\begin{lemma}
There holds the following identity
$$\int_{-\infty}^{\infty} \cos(m\pi y)\frac{1}{\sqrt{4k\pi t}}e^{-\frac{(x-y)^2}{4kt}} dy = \cos(m\pi x)e^{-m^2\pi^2kt}.$$
\end{lemma}
\textbf{Proof.} First, by Taylor expansion of $\cos(m\pi y)$ we have
\begin{eqnarray*}
\int_{-\infty}^{\infty} \cos(m\pi y)\frac{1}{\sqrt{4k\pi t}}e^{-\frac{y^2}{4kt}} dy
& = & \int_{-\infty}^{\infty} \sum_{i=0}^{\infty}(-1)^i \frac{(m\pi y)^{2i}}{(2i)!}
\frac{1}{\sqrt{4k\pi t}}e^{-\frac{y^2}{4kt}} dy \\
& = & \sum_{i=0}^{\infty}\frac{(-1)^i (m\pi)^{2i}}{(2i)!} \int_{-\infty}^{\infty} y^{2i}
\frac{1}{\sqrt{4k\pi t}}e^{-\frac{y^2}{4kt}} dy \\
& = & \sum_{i=0}^{\infty}\frac{(-1)^i (m\pi)^{2i}}{(2i)!} (4kt)^i\prod_{j=1}^i \frac{2j-1}{2} \\
& = & \sum_{i=0}^{\infty}\frac{(-1)^i (m\pi)^{2i}}{i!}(kt)^i= e^{-m^2\pi^2kt}.
\end{eqnarray*}
Thus we have
\begin{eqnarray*}
& & \int_{-\infty}^{\infty} \cos(m\pi y)\frac{1}{\sqrt{4k\pi t}}e^{-\frac{(x-y)^2}{4kt}} dy \\
 & =& \int_{-\infty}^{\infty} \cos(m\pi (x-y))\frac{1}{\sqrt{4k\pi t}}e^{-\frac{y^2}{4kt}} dy \\
& = &  \int_{-\infty}^{\infty} \left[\cos(m\pi x) \cos(m\pi y)+\sin(m\pi x) \sin(m\pi y)\right]
\frac{1}{\sqrt{4k\pi t}}e^{-\frac{y^2}{4kt}} dy\\
& = & \cos(m\pi x) \int_{-\infty}^{\infty} \cos(m\pi y)\frac{1}{\sqrt{4k\pi t}}e^{-\frac{y^2}{4kt}} dy\\
& = & \cos(m\pi x) e^{-m^2\pi^2kt}.
\end{eqnarray*}
\hspace{\fill} $\Box$\\
By Lemma 1 and (\ref{eq:301}) one can easily get another solution form of (\ref{eq:101})
\begin{equation}
\label{eq:302}
u(x,t)= \sum_{m=0}^{\infty}\left[e^{-m^2\pi^2kt}a_m \cos(m\pi x) +b_m\cos(m\pi x)\int_0^t h(s)e^{-m^2\pi^2k(t-s)}ds  \right],
\end{equation}
where for $m=0$ there hold $a_0=\int_0^1 \mu_0(x)dx$, $b_0=\int_0^1 f(x)dx$, and for $m>0$ there hold
$a_m= 2\int_0^1 \mu_0(x) \cos(m\pi x) dx$,  $b_m=2\int_0^1 f(x) \cos(m\pi x) dx$.
By using the final data $\mu_T(x)=u(x,T)$ we have
$$\mu_T(x)= \sum_{m=0}^{\infty}\left[e^{-m^2\pi^2kT}a_n \cos(m\pi x) +b_n\cos(m\pi x)\int_0^T h(s)e^{-m^2\pi^2k(T-s)}ds  \right].$$
Denoting $c_0=\int_0^1 \mu_T(x)dx$ and $c_n=2\int_0^1 \mu_T(x) \cos(m\pi x) dx, m>0$ and integrating both sides of
the above equation with $\cos(m\pi x)$ we obtain
\begin{equation}
\label{eq:303}
b_m = \frac{e^{m^2\pi^2kT}c_m- a_m}{\int_0^T h(s)e^{m^2\pi^2 k s }ds}.
\end{equation}
We can use the singular decomposition to solve $b_m$. In fact, if we define linear operator $K$ as
$$Kf(x):=\sum_{n=0}^{\infty}\frac{\int_0^T h(s)e^{-m^2\pi^2k(T-s)}ds}{\int_0^T h(x) ds}b_m\cos(m\pi x)$$
then we have the singular values (or eigenvalues) $\{\sigma_m\}$ with
$\sigma_m=\frac{\int_0^T h(s)e^{-m^2\pi^2k(T-s)}ds}{\int_0^T h(x) ds}$
and corresponding eigenvector $\{\cos(m\pi x)\}$. Since $|\sigma_m|\leq 1$ and $\lim_{m\rightarrow \infty} \sigma_m=0$.
The problem is a general linear operator equation in inverse problem. Lots of regularization methods, such as Tikhonov regularization,
Landerweber iteration, etc., can be used to solve this problem. In this paper we shall discuss about the frequency
cut-off method to solve the problem, references for other methods can be found in \cite{EHN:1996,AK:1996}.  In fact,
we only need to solve $b_m$ in (\ref{eq:303}). We denote $c_0^{\delta}=\int_0^1 \mu_T^{\delta}(x)dx$,
$c_m^{\delta}=\int_0^1 \mu_T^{\delta}(x) \cos(m\pi x) dx, m>0$, and
\begin{equation}
\label{eq:304} \|\mu_T^{\delta}(\cdot)-\mu_T(\cdot)\| = \left(\int_0^1 \left[\sum_{m=0}^{\infty}
(c_m^{\delta}-c_m)\cos(n\pi x) \right]^2dx\right)^{\frac{1}{2}} \leq \delta.
\end{equation}
We introduce the frequency cut-off method to solve $b_m^{\delta}$ by
\begin{equation}
\label{eq:305}
b_m^{\delta} = \frac{\chi_{\vartheta}(e^{m^2\pi^2kT}c_m^{\delta}- a_m)}{\int_0^T h(s)e^{m^2 \pi^2 k s }ds},
\end{equation}
where $\chi_{\vartheta}$ is the discrete version of characteristic function defined in section 2, that is
$\chi_{\vartheta}=1$ for $m\leq \vartheta \in N$ and $\chi_{\vartheta}=0$ for $m>\vartheta$. Thus
we can reconstruct $f(x)$ with $f^{\delta}(x)$
$$f^{\delta}(x)=\sum_{m=0}^{\vartheta} b_m^{\delta} cos (m \pi x)$$
and
$$u^{\delta}(x,t)= \sum_{m=0}^{\infty}e^{-m^2\pi^2kt}a_m cos(m\pi x)
+\sum_{m=0}^{\vartheta}b^{\delta}_mcos(m\pi x)\int_0^t h(s)e^{-m^2\pi^2k(t-s)}ds.$$
Based
on (\ref{eq:305}) we have the convergence theorem as follows
\begin{theorem}
Let $u(x,t)$ be the exact temperature history of (\ref{eq:101}), $h(t)\not\equiv 0$ is identically
nonpositive or nonnegative in $[0,T]$ and $\mu_T^{\delta}(x)$ be the measured final temperature
satisfying (\ref{eq:304}). $f(x)$ satisfies $\|f\|_{H^p(0,1)}\leq M$, where $\|f\|_{H^p(0,1)}$ is defined as
$$\|f\|_{H^p(0,1)}= \left(\sum_{m=0}^{\infty} (1+m^2)^p b_m^2\right)^{\frac{1}{2}}.$$
Let $b_m^{\delta}$ defined by (\ref{eq:305}) and
$\vartheta \sim \lfloor \sqrt{\frac{1}{(1+\sigma)kT}
\left[ln\frac{M}{\delta}(ln\frac{M}{\delta})^{-\frac{1+\sigma}{2}p}\right]}\rfloor$, $\sigma\geq 0$
then there holds
\begin{equation}
\label{eq:th201} {\| f^{\delta }-f
\|}^{2} \leq C(ln\frac{M}{\delta})^{-p}
\left(  M^{\frac{2}{1+\sigma}}\delta^{\frac{2\sigma}{1+\sigma}}+ M^{2}(\frac{ln\frac{M}{\delta}}
{ln\frac{M}{\delta}(ln\frac{M}{\delta})^{-\frac{1+\sigma}{2}p}})^{p} \right),
\end{equation}
and
\begin{equation}
\label{eq:th202}
{\| u^{\delta }(\cdot,t)-u(\cdot,t)\|}^{2}
\leq C_h(t)(ln\frac{M}{\delta})^{-p}
\left(  M^{\frac{2}{1+\sigma}}\delta^{\frac{2\sigma}{1+\sigma}}+ M^{2}(\frac{ln\frac{M}{\delta}}
{ln\frac{M}{\delta}(ln\frac{M}{\delta})^{-\frac{1+\sigma}{2}p}})^{p} \right),
\end{equation}
for $\delta \rightarrow 0$,
where $C$ is a constant independent of $\delta$ and $M$ and $C_h(t)=C(\int_0^t |h(s)|ds)^2$.
\end{theorem}
\textbf{Proof.}
\begin{eqnarray*}
\|f^{\delta}-f\|^2 & = & \|\sum_{m=0}^{\vartheta} (b_m^{\delta}-b_m) cos (m \pi \cdot)-
\sum_{m=\vartheta+1}^{\infty}b_m cos (m \pi \cdot)\|^2 \\
& \leq & 2 \|\sum_{m=0}^{\vartheta} (b_m^{\delta}-b_m) cos (m \pi \cdot)\|^2
+ 2\|\sum_{m=\vartheta+1}^{\infty}b_m cos (m \pi \cdot)\|^2 \\
& = & 2\|\sum_{m=0}^{\vartheta}
\frac{e^{m^2\pi^2kT}(c_m^{\delta}- c_m)}{\int_0^T h(s)e^{m^2 \pi^2 k s }ds} cos (m \pi \cdot)\|^2
+ 2\|\sum_{m=\vartheta+1}^{\infty}b_m cos (m \pi \cdot)\|^2 \\
& \leq & \frac{2}{C_h^2}e^{\vartheta^2\pi^2kT} \delta^2 + 2 \vartheta^{-2p}M^2 \\
& \leq & 2(ln\frac{M}{\delta})^{-p}
\left( \frac{1}{C_h^2} M^{\frac{2}{1+\sigma}} \delta^{\frac{2\sigma}{1+\sigma}}+ CM^{2}(\frac{ln\frac{M}{\delta}}
{ln\frac{M}{\delta}(ln\frac{M}{\delta})^{-\frac{1+\sigma}{2}p}})^{p} \right).
\end{eqnarray*}
\begin{eqnarray*}
{\| u^{\delta }(\cdot,t)-u(\cdot,t)
\|}^{2} & = & \|(f^{\delta }-f)\int_0^t h(s)e^{-m^2\pi^2k(t-s)}ds\|^2 \\
& \leq & 2(\int_0^t|h(s)|ds)^2(ln\frac{M}{\delta})^{-p}
\left(  M^{\frac{2}{1+\sigma}} \delta^{\frac{2\sigma}{1+\sigma}}+ M^{2}(\frac{ln\frac{M}{\delta}}
{ln\frac{M}{\delta}(ln\frac{M}{\delta})^{-\frac{1+\sigma}{2}p}})^{p} \right),
\end{eqnarray*}
which proves the theorem. \hspace{\fill} $\Box$
\begin{remark}
We see from Theorem 3, that the convergence results are similar to the results in Section 2. Actually, we
can also design the similar iterative method like we design in Section 2. And we can introduce the discrepancy principle
as the \emph{a posteriori} stopping rule. If (\ref{eq:301}) is used for reconstruction then the corresponding functions are
extended to periodical functions. Thus in the numerical calculation, by using iterative method in section 2
the fast Fourier transform (FFT) can be considered for the periodical function in $\RR$. 
\end{remark}
\section{Numerical experiment}
In this section, we present some numerical experiments on reconstruction of the source term with the final
measurement $\mu_T^{\delta}(x)$ for $T=1$. We separate the span $[0,1]$ for $x$ variable into an equidistance grid
$0=x_0<\cdots< x_i < \cdots<x_{N_1}=1$ ($x_i=ih, h=0.02, N_1=50$),
and the span $[0,1]$ for $t$ variable into an equidistance grid
$0=t_0<\cdots<t_j<\cdots<t_{N_2}=1$ ($t_j=jl, l=0.05, N_2=20$).
We produce the random noise as follows
$$\mu_T^{\delta}(x_i)=\mu_T^{\delta}(x_i)+ 2(\mbox{rand}(0,1)-0.5)*noiselv*\mu_T^{\delta}(x_i),$$
where rand(0,1) denotes the uniformly distributed pseudo-random numbers in [0,1] generated by Matlab software
and $noiselv$ is a positive number between 0 and 1 for noise level. The noise $\delta$ is calculated by
numerical calculation of $L^2(0,1)$ norm (by first approximating the function with spline interpolation and then
using the integral algorithm). \\
We only consider the numerical implementation of reconstruction of the source term in problem (\ref{eq:101}) although
numerical method for \eqnref{eq:103} can be similarly implemented.
\subsection{Example 1} 
\begin{figure}
   \begin{center}
         \includegraphics[width=2.8in]{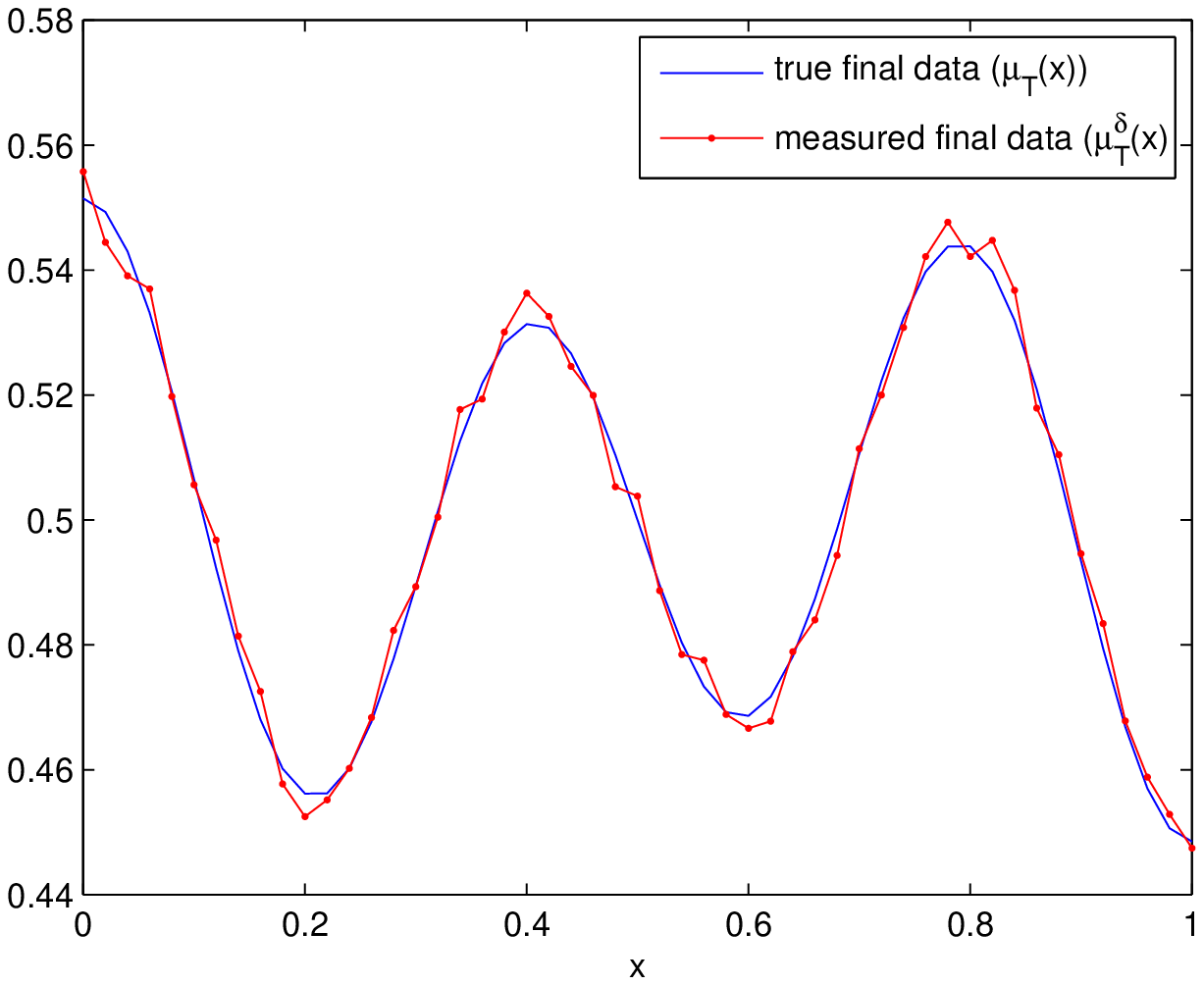}
         \includegraphics[width=2.8in]{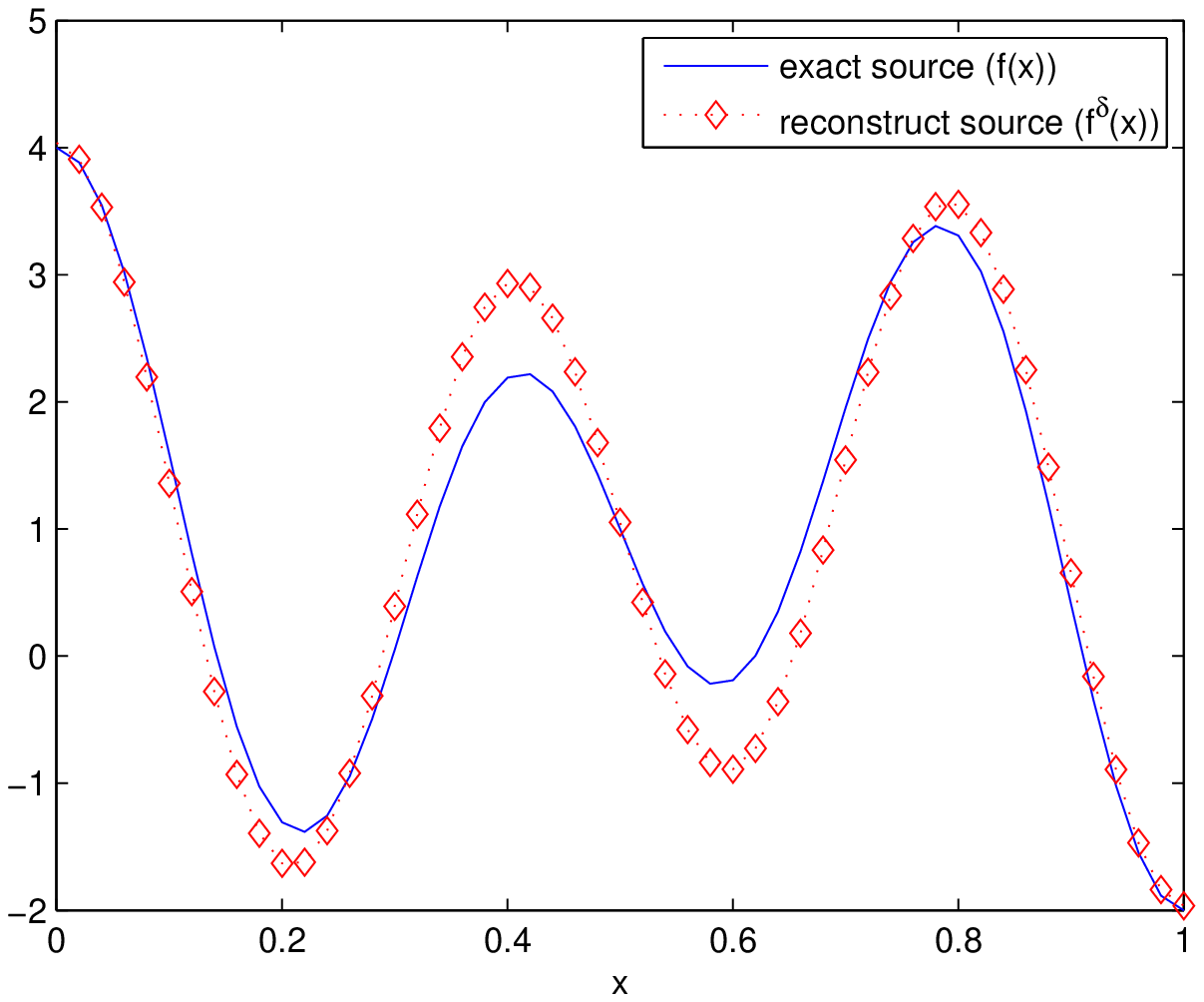}
         \includegraphics[width=2.8in]{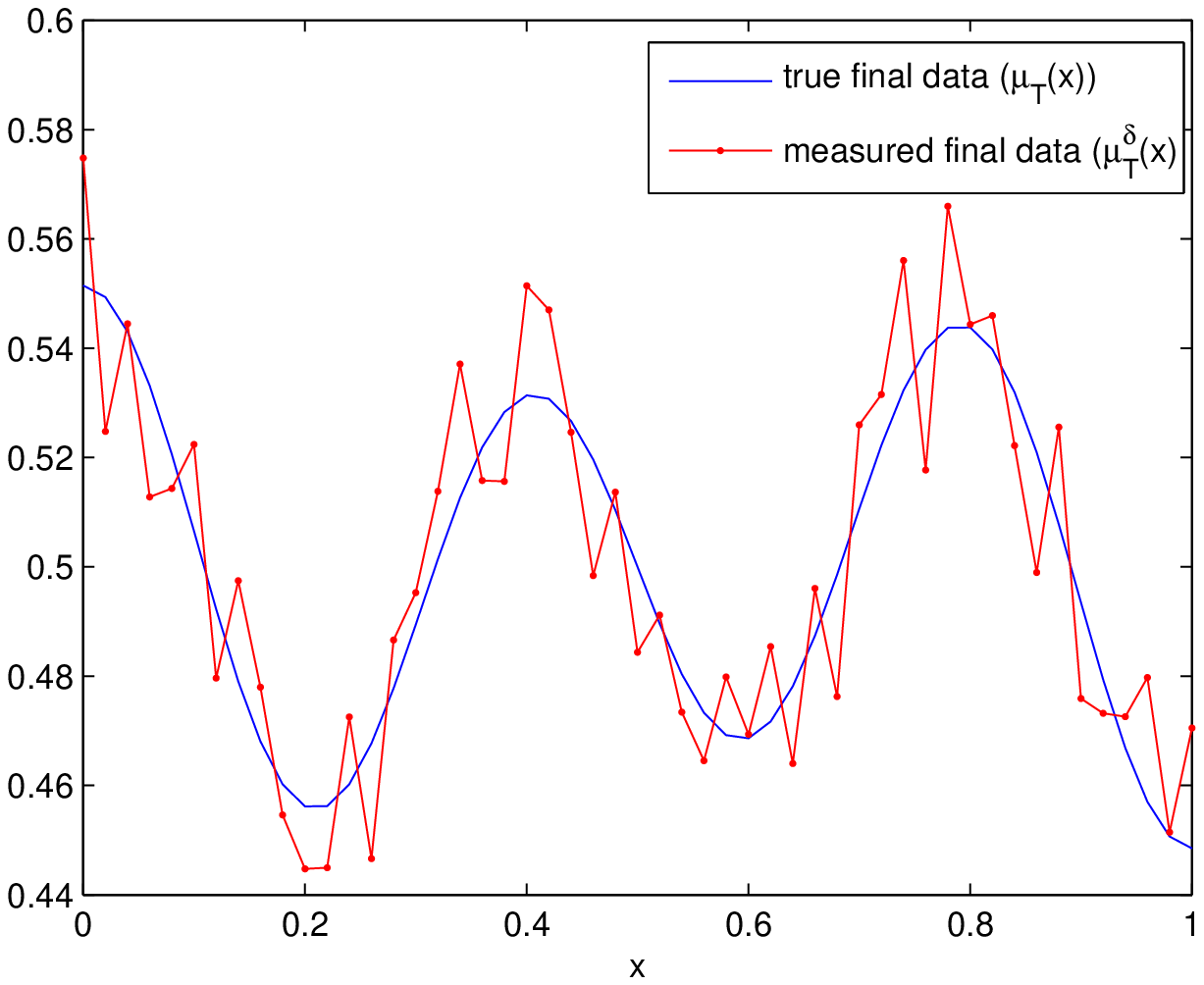}
         \includegraphics[width=2.8in]{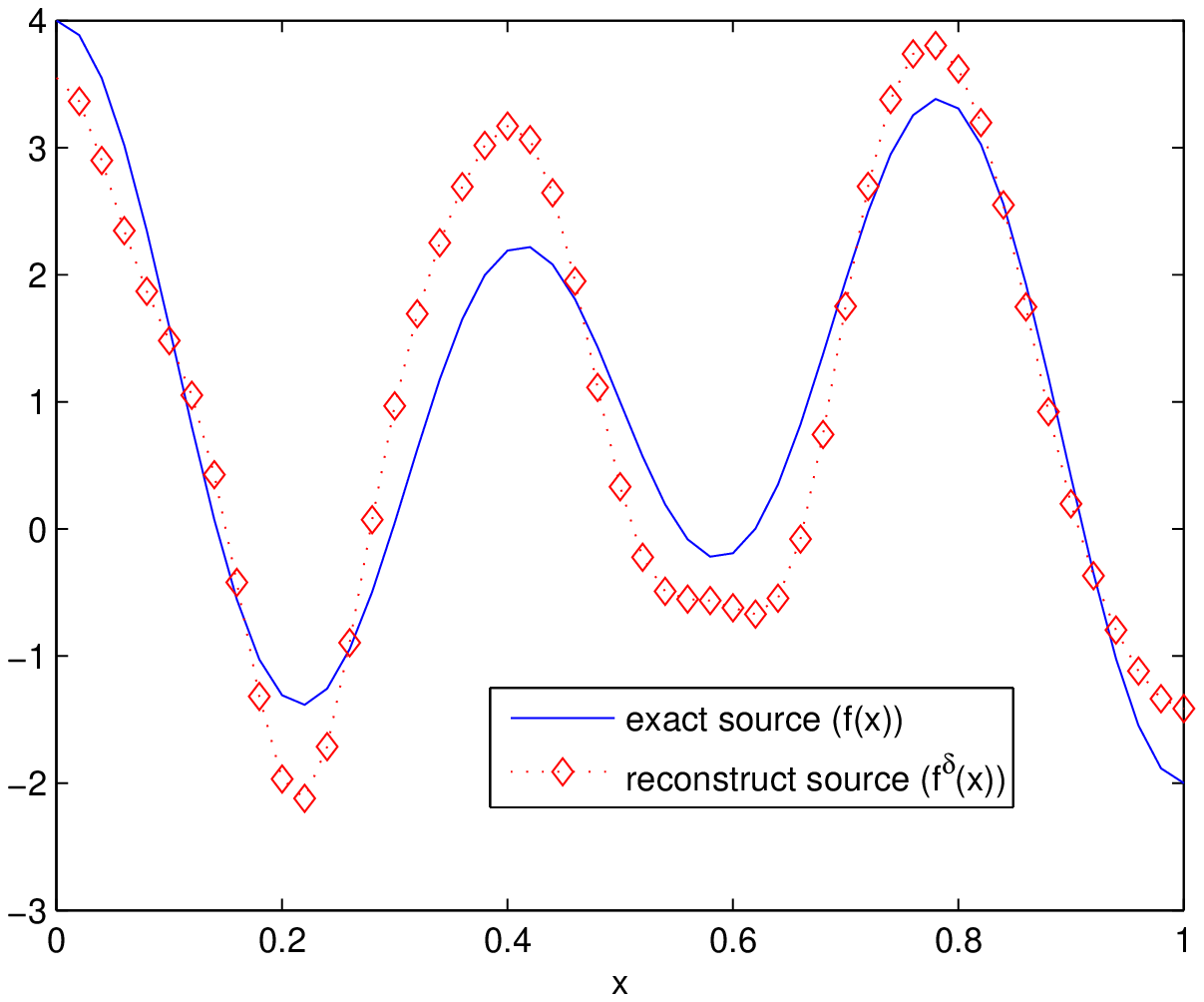}
   \end{center}
    \caption{Comparison of convergence results. The top two with the noise level $1\%$ and $\vartheta=12$, while the
        bottom two with the noise level $5\%$ and $\vartheta=18$. The left are the final data and corresponding
        measurement error data. The other parameters are $p=1/3$, $\sigma=0.2$.
       }\label{fig:1}
\end{figure}

\begin{figure}
   \begin{center}
         \includegraphics[width=2.8in]{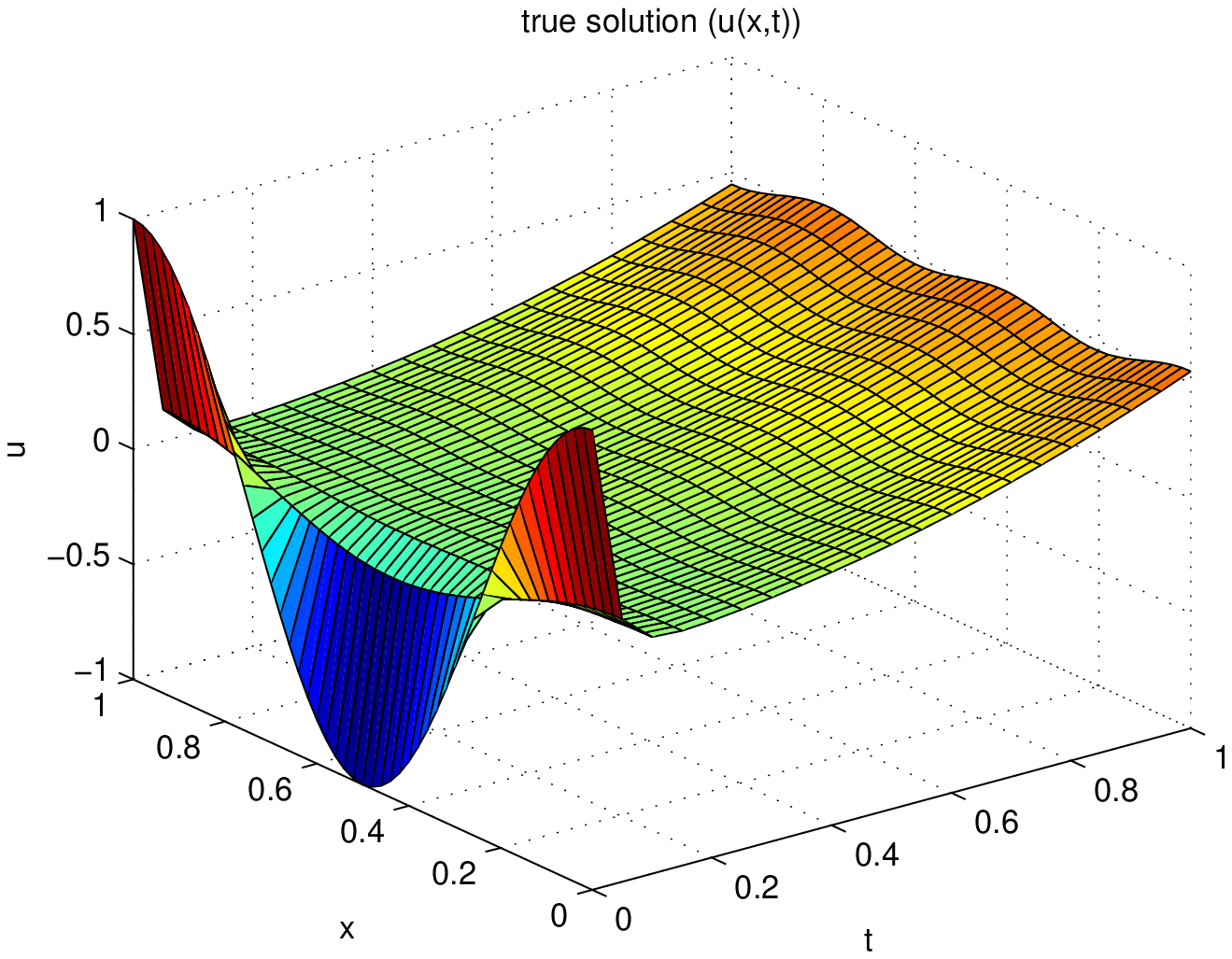}
         \includegraphics[width=2.8in]{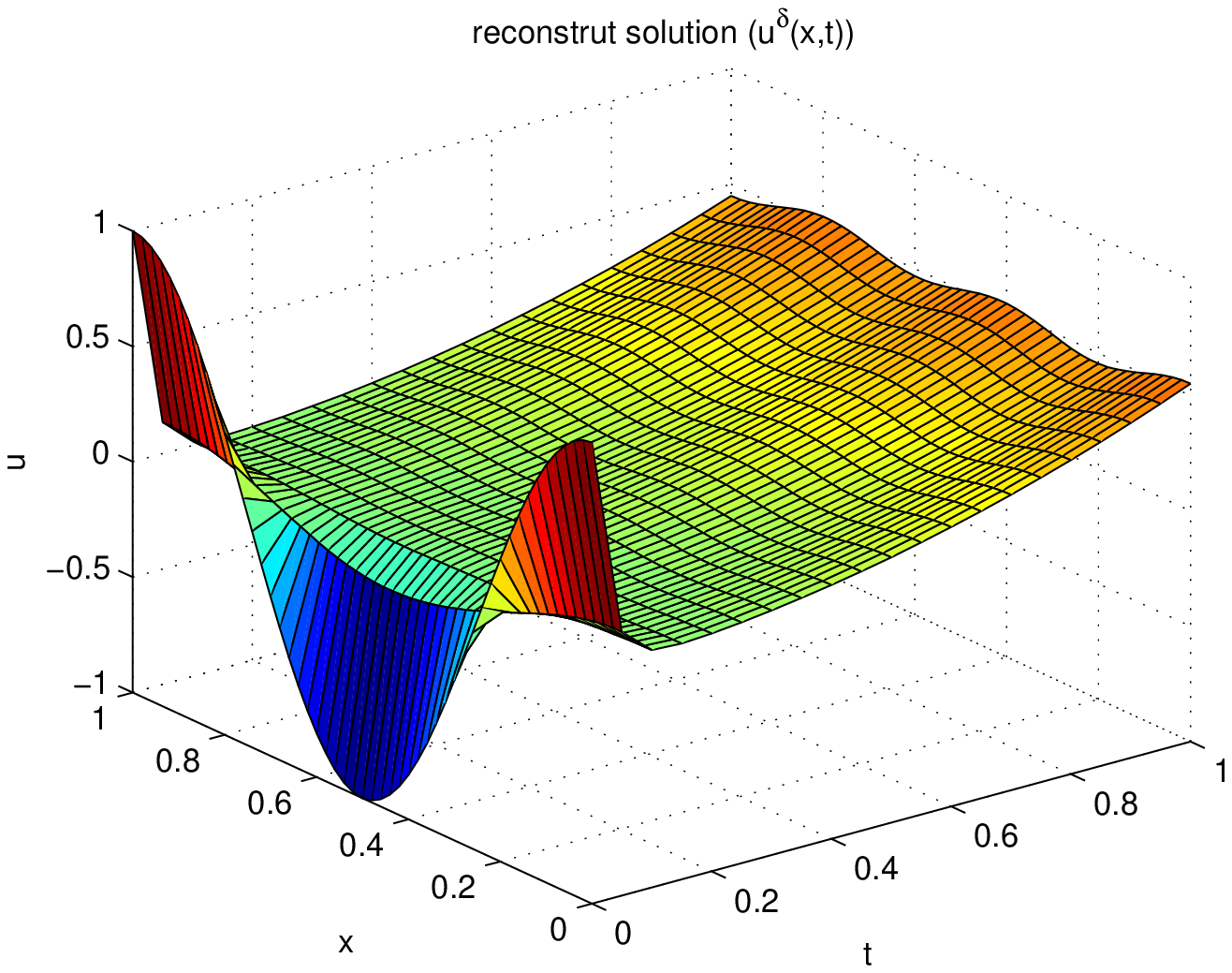}
   \end{center}
    \caption{Comparison of true solution and reconstruct solution. The noise level is $1\%$ and $\vartheta=12$, $p=1/3$, $\sigma=0.2$.
       }\label{fig:2}
\end{figure}

\begin{figure}
   \begin{center}
       \includegraphics[width=2.8in]{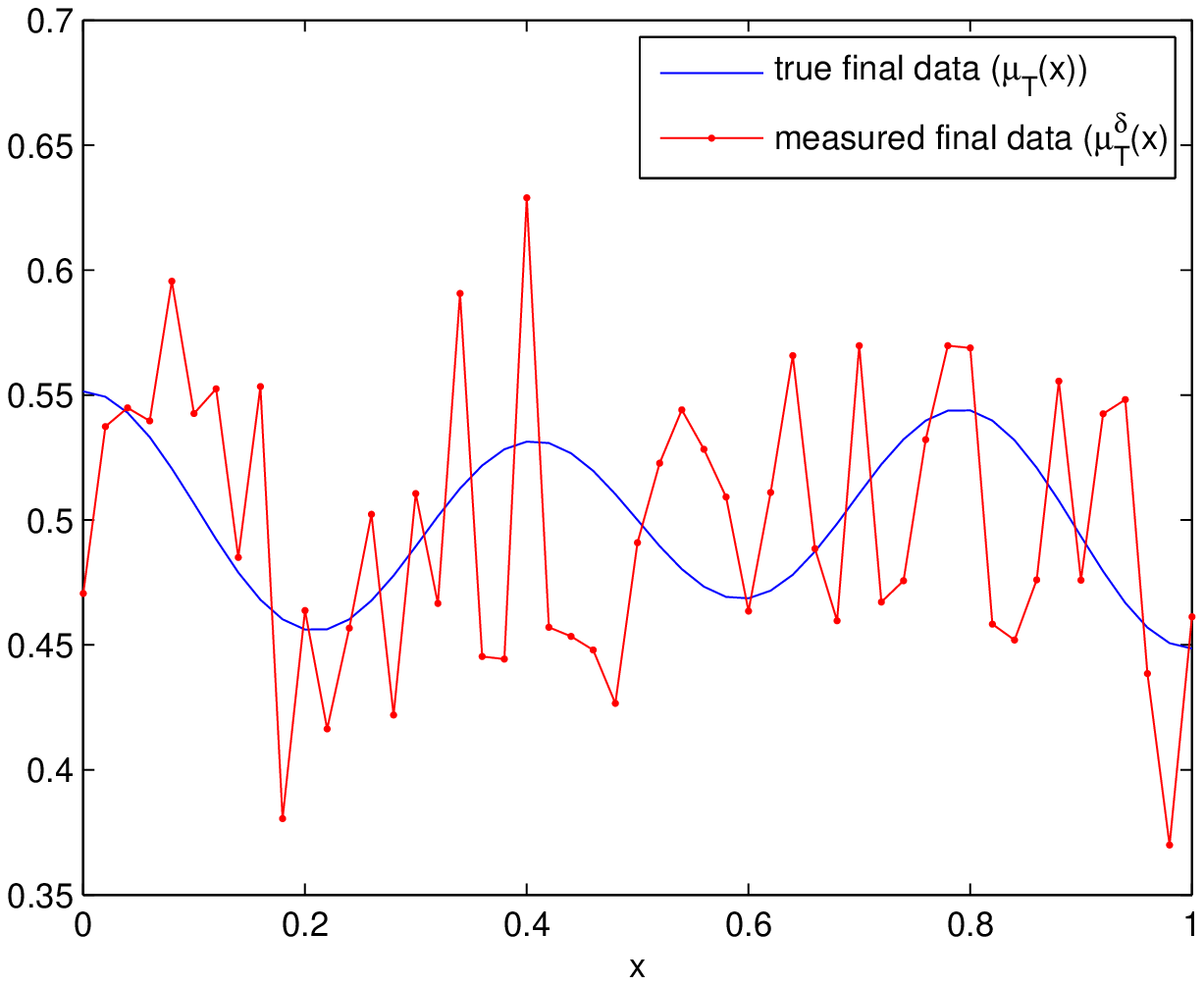}
       \includegraphics[width=2.8in]{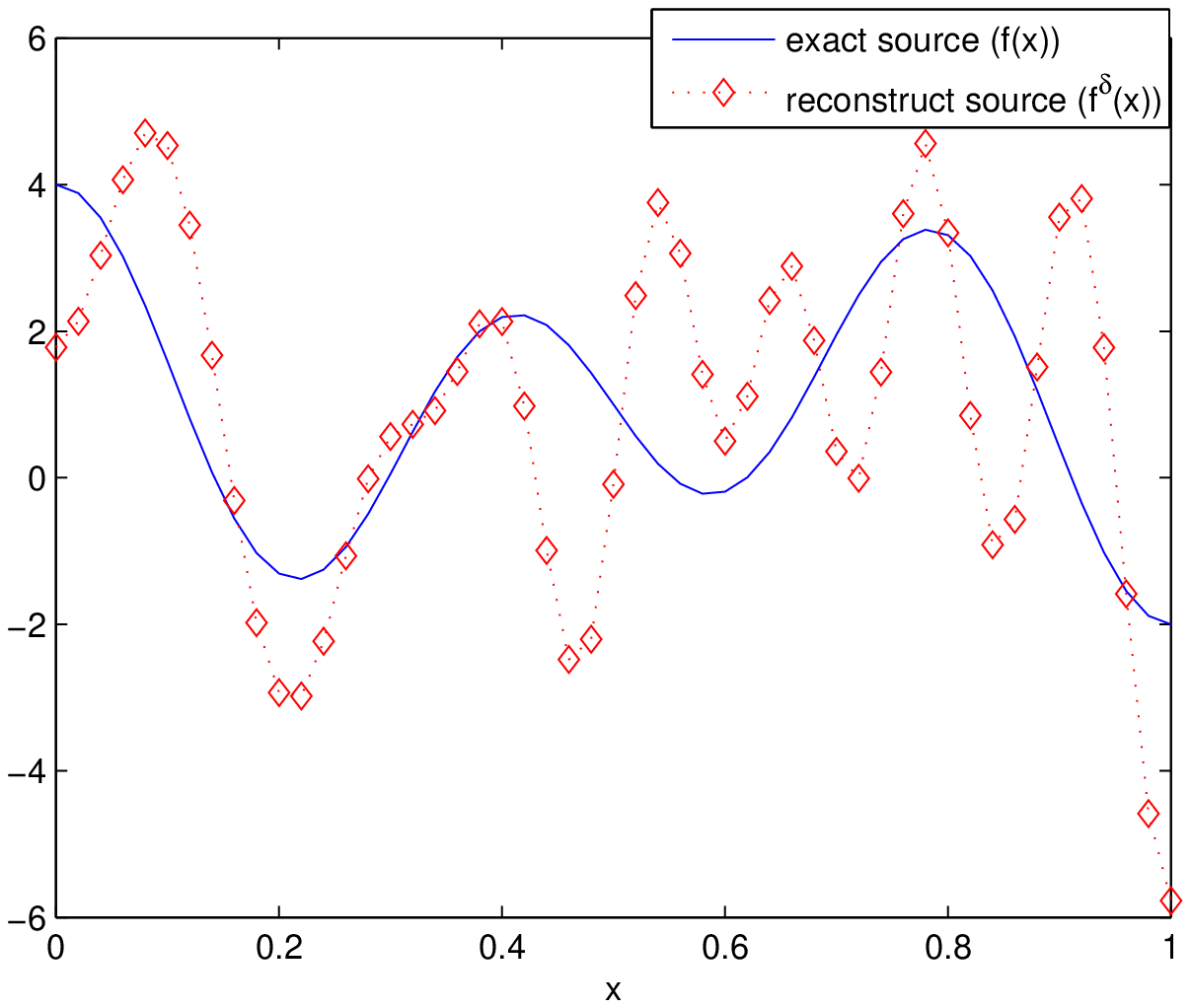}
        \includegraphics[width=2.8in]{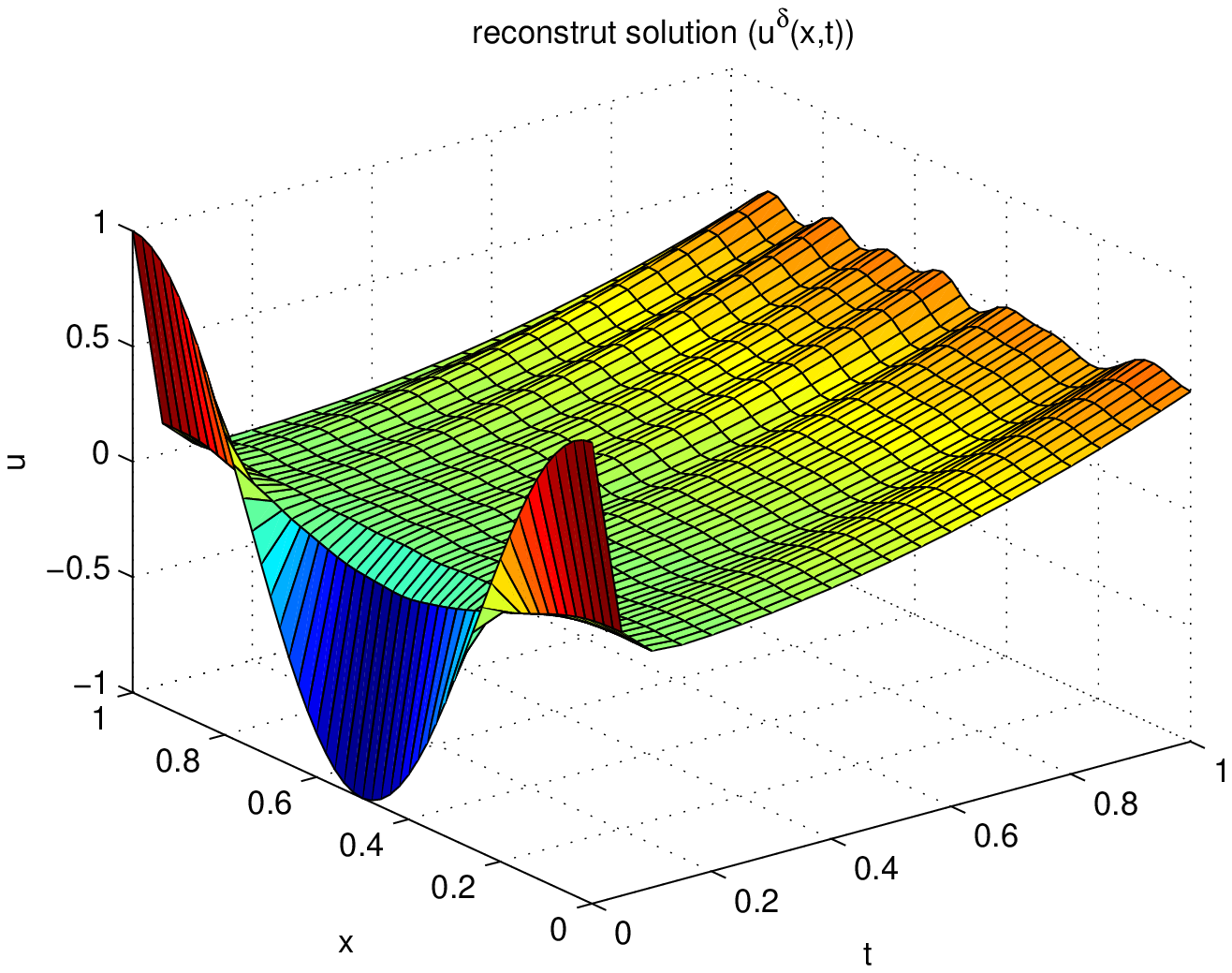}
   \end{center}
    \caption{Convergence results with noise level $20\%$ and $\vartheta=18$, $p=1/3$, $\sigma=0.2$.
       }\label{fig:3}
\end{figure}

Set $f(x)=1+\cos(3\pi x) + 2\cos(5\pi x)$, $\mu_0(x)=\cos(2\pi x)$, $h(t)=t$ and $k=1$ then
the solution of (\ref{eq:101}) is
$$u(x,t)=\cos(2\pi x)e^{-4\pi^2t} + \frac{t^2}{2} + \frac{9\pi^2 t-1+e^{-9\pi^2t}}{81\pi^4}\cos(3\pi x)
+ 2\frac{25\pi^2t-1+e^{-25\pi^2t}}{125\pi^4} \cos(5\pi x),$$
and the final measurement at $T=1$ is
$$\mu_T(x,t)=\frac{1}{2} +\cos(2\pi x)e^{-4\pi^2} +  \frac{9\pi^2 -1+e^{-9\pi^2}}{81\pi^4}\cos(3\pi x)
+ 2\frac{25\pi^2-1+e^{-25\pi^2}}{125\pi^4} \cos(5\pi x).$$
We choose $p=1/3$ and $\sigma=0.2$ in regularization method (\ref{eq:305}). For the numerical calculation of the first
term in (\ref{eq:302}), i.e.,
$\sum_{m=0}^{\infty}e^{-m^2\pi^2t}\cos(m\pi x) \int_0^1 \mu_0(y) \cos(m\pi y) dy$, 
 we choose a sufficiently large number $m<3\vartheta$ to numerically approximate it. Thus the numerical implementation of $u^{\delta}(x,t)$ is
$$u^{\delta}(x,t)= \sum_{m=0}^{3\vartheta}e^{-m^2\pi^2t}a_m cos(m\pi x)
+\sum_{m=0}^{\vartheta}b^{\delta}_mcos(m\pi x)\int_0^t h(s)e^{-m^2\pi^2(t-s)}ds.$$
Table 1 shows the numerical results for different choice of $\vartheta$ and error level. We see from the
table that the reconstruction of the solution $u^{\delta}$ has more accuracy than the source $f(x)$. Fig. \ref{fig:1}
shows the performance of reconstruction of source term $f(x)$ under different final measurements $\mu_T^{\delta}(x)$
while Fig. \ref{fig:2} gives the comparison between the true solution and the numerical solution with noise level
$1\%$. It is clearly in this figure that the solution is not affected that much compared with the source term 
under the measurement noise of $\mu_T(x)$.
Fig. \ref{fig:3} shows the numerical results under the noise level $20\%$. Since the noise level is quite
high, the numerical method can not produce good approximation solution.
\begin{table}[h]
\label{tab:401} \caption{Convergence results. The parameter $p=1/3$, $\sigma=0.2$, $M=1.870888$.}
\begin{center}
\begin{tabular}{lcc}
\hline\noalign{\smallskip}
 & $noise \ \ level (1\%)$&$noise \ \ level (5\%)$ \\
$\delta$ & 0.003035 & 0.013690 \\
$\vartheta=$ &   12 | 24 | 36 & 6 | 12 | 18 \\
\noalign{\smallskip} \hline \noalign{\smallskip}
$\|f^{\delta}-f\|$                     & 0.383497 0.423780 0.673838& 0.346827 0.617770 0.580466\\
$\|u^{\delta}(\cdot,t)-u(\cdot,t)\|$   & 0.017470 0.044817 0.183869& 0.078773 0.143316 0.133631\\
\noalign{\smallskip} \hline \noalign{\smallskip}
\noalign{\smallskip} \hline \noalign{\smallskip}
&$noise \ \ level (10\%)$ & $noise \ \ level (20\%)$ \\
$\delta$ &  0.027804 & 0.055623 \\
$\vartheta=$ &  6 | 12 | 18 & 6 | 12 | 18 \\
\noalign{\smallskip} \hline \noalign{\smallskip}
$\|f^{\delta}-f\|$                     & 0.560617 0.510970 1.120361& 0.744648 1.542746 2.033309  \\
$\|u^{\delta}(\cdot,t)-u(\cdot,t)\|$   & 0.177541 0.187410 0.347001& 0.415692 0.564468 0.727249  \\
\noalign{\smallskip} \hline \noalign{\smallskip}
\end{tabular}
\end{center}
\end{table}

\clearpage{}
\subsection{Example 2}
\begin{figure}[h]
   \begin{center}
         \includegraphics[width=2.8in]{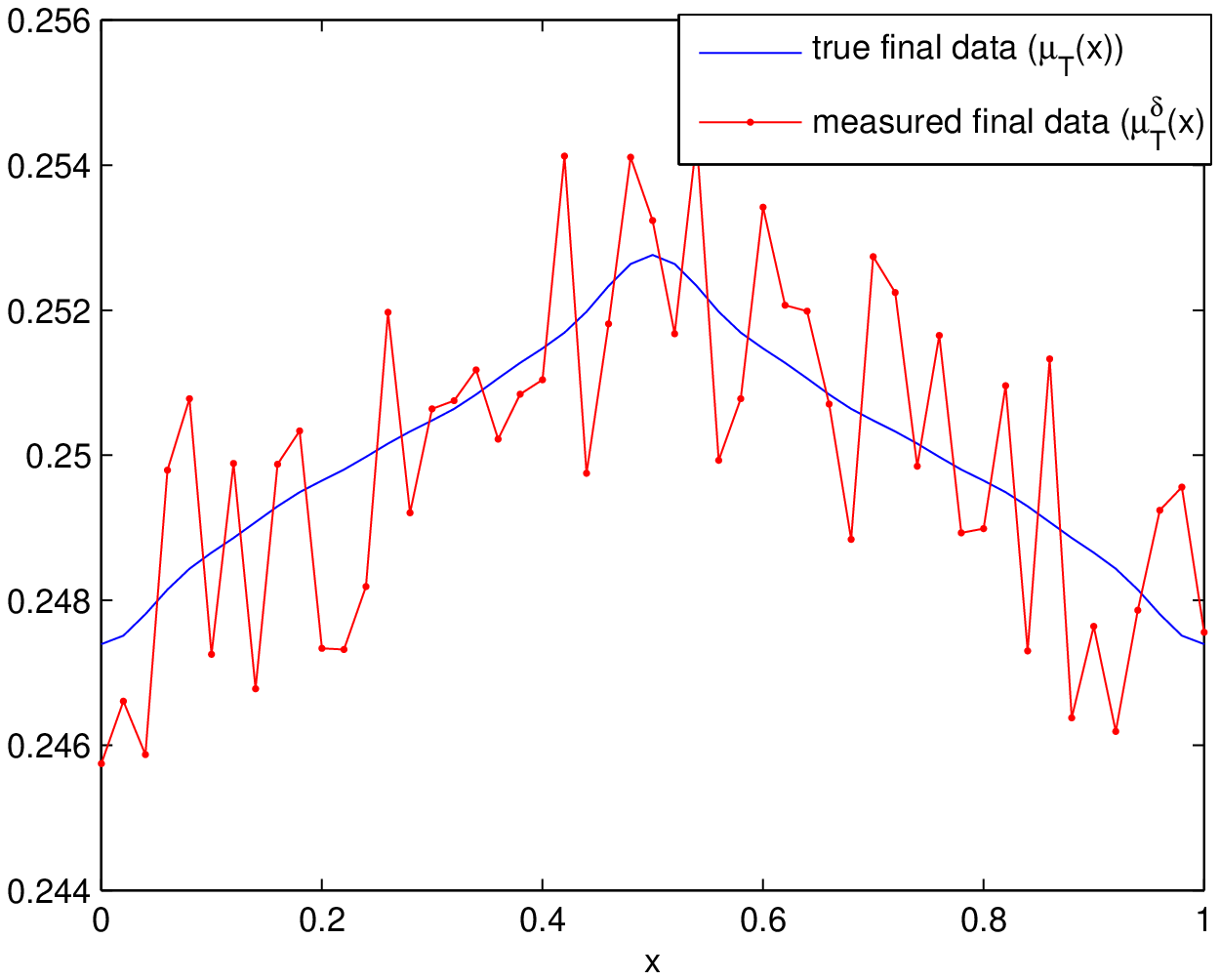}
         \includegraphics[width=2.8in]{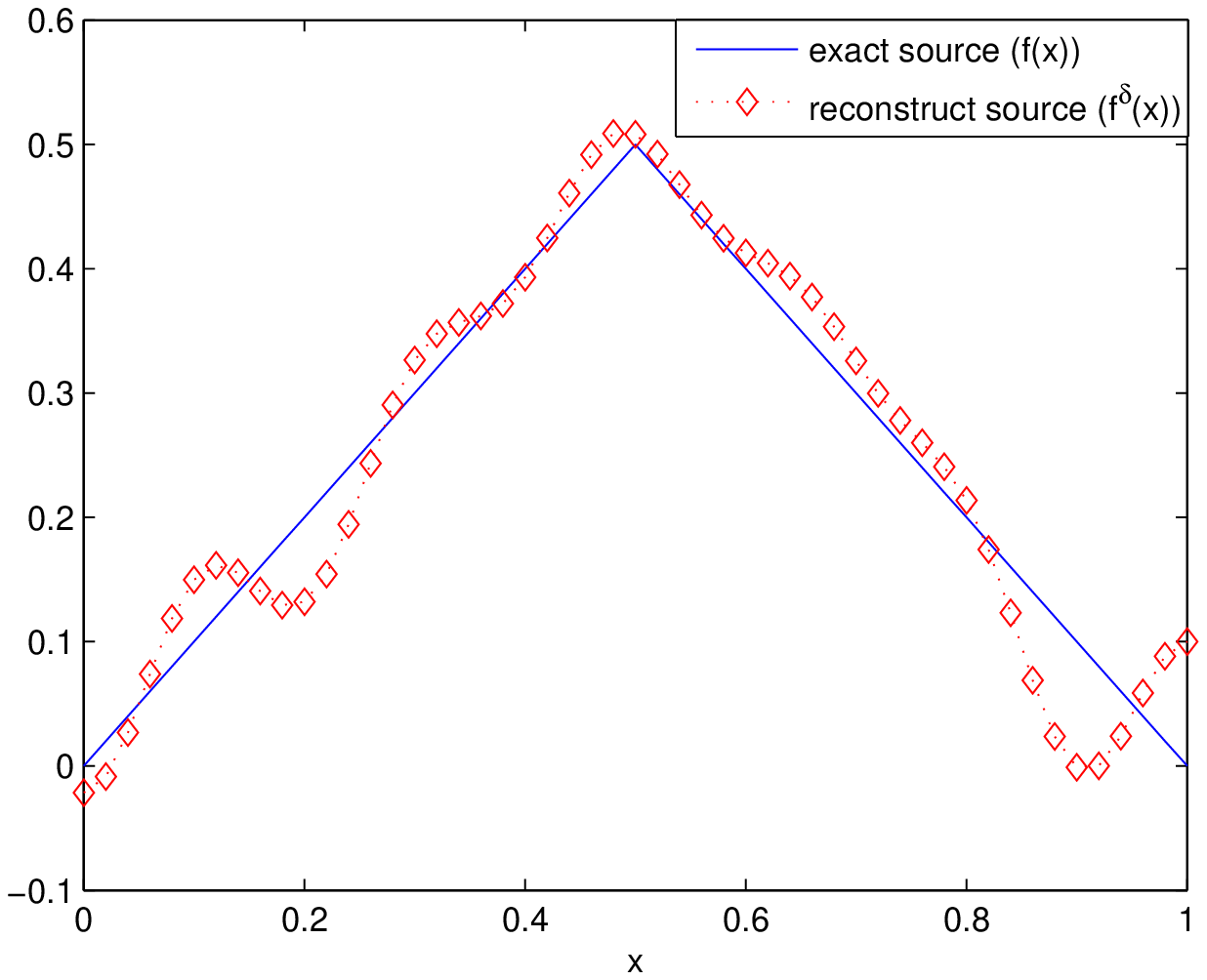}
         \includegraphics[width=2.8in]{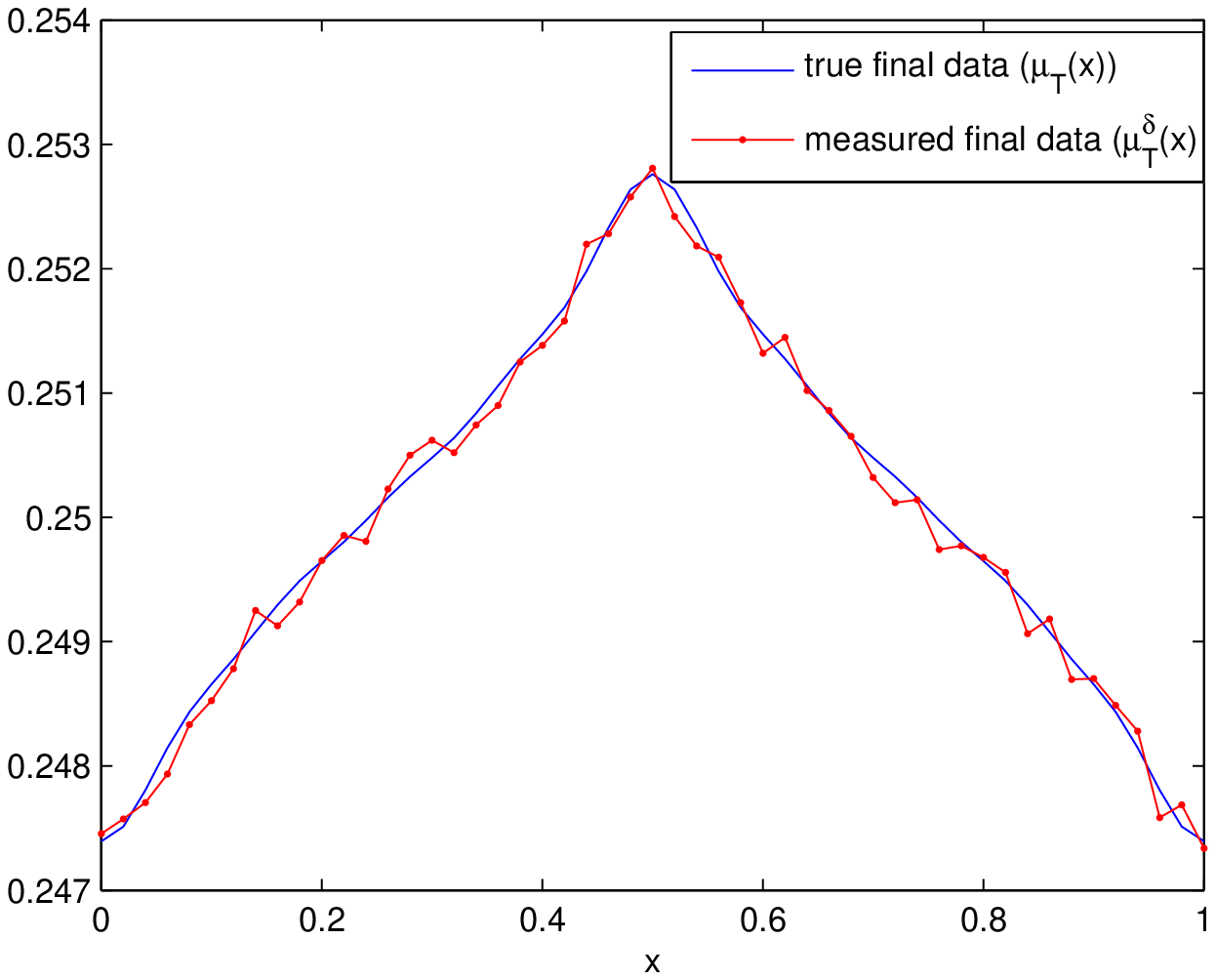}
         \includegraphics[width=2.8in]{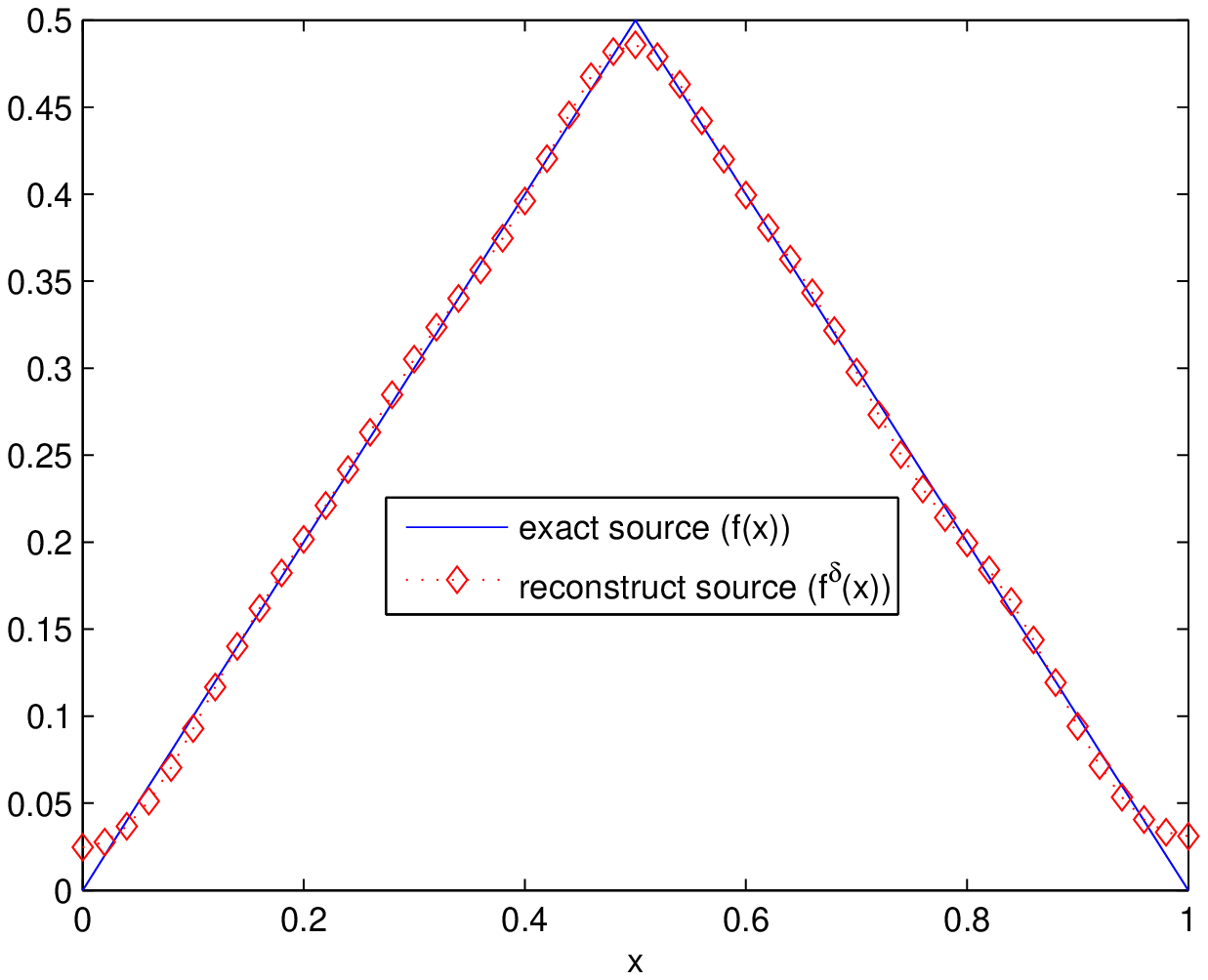}
   \end{center}
    \caption{Comparison of convergence results. The top two with the noise level $1\%$, while the
        bottom two with the noise level $0.1\%$. The left are the final data and corresponding
        measurement error data. The other parameters are $\vartheta=12$, $p=1/3$, $\sigma=0.2$.
       }\label{fig:4}
\end{figure}

\begin{figure}[h]
   \begin{center}
         \includegraphics[width=2.8in]{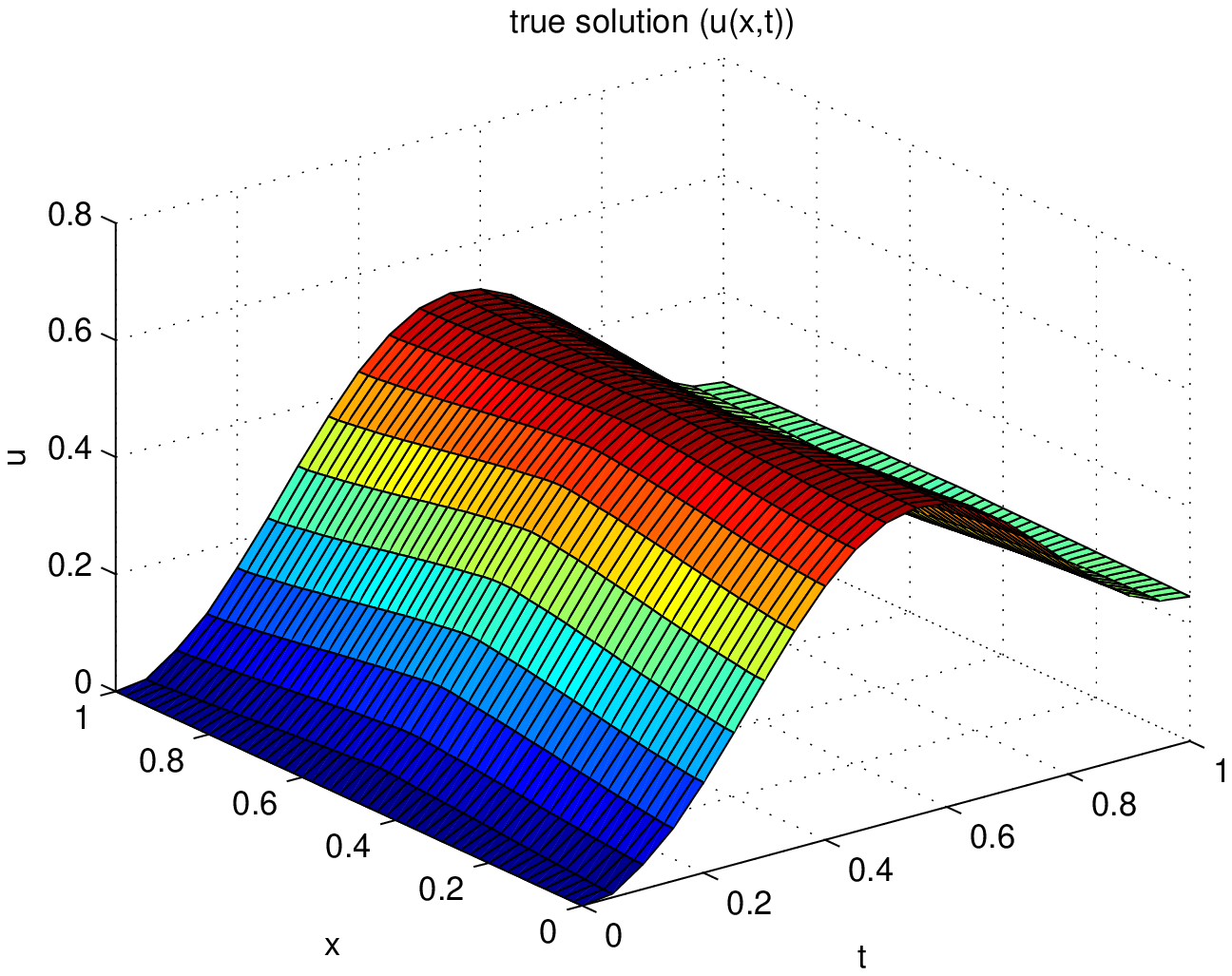}
         \includegraphics[width=2.8in]{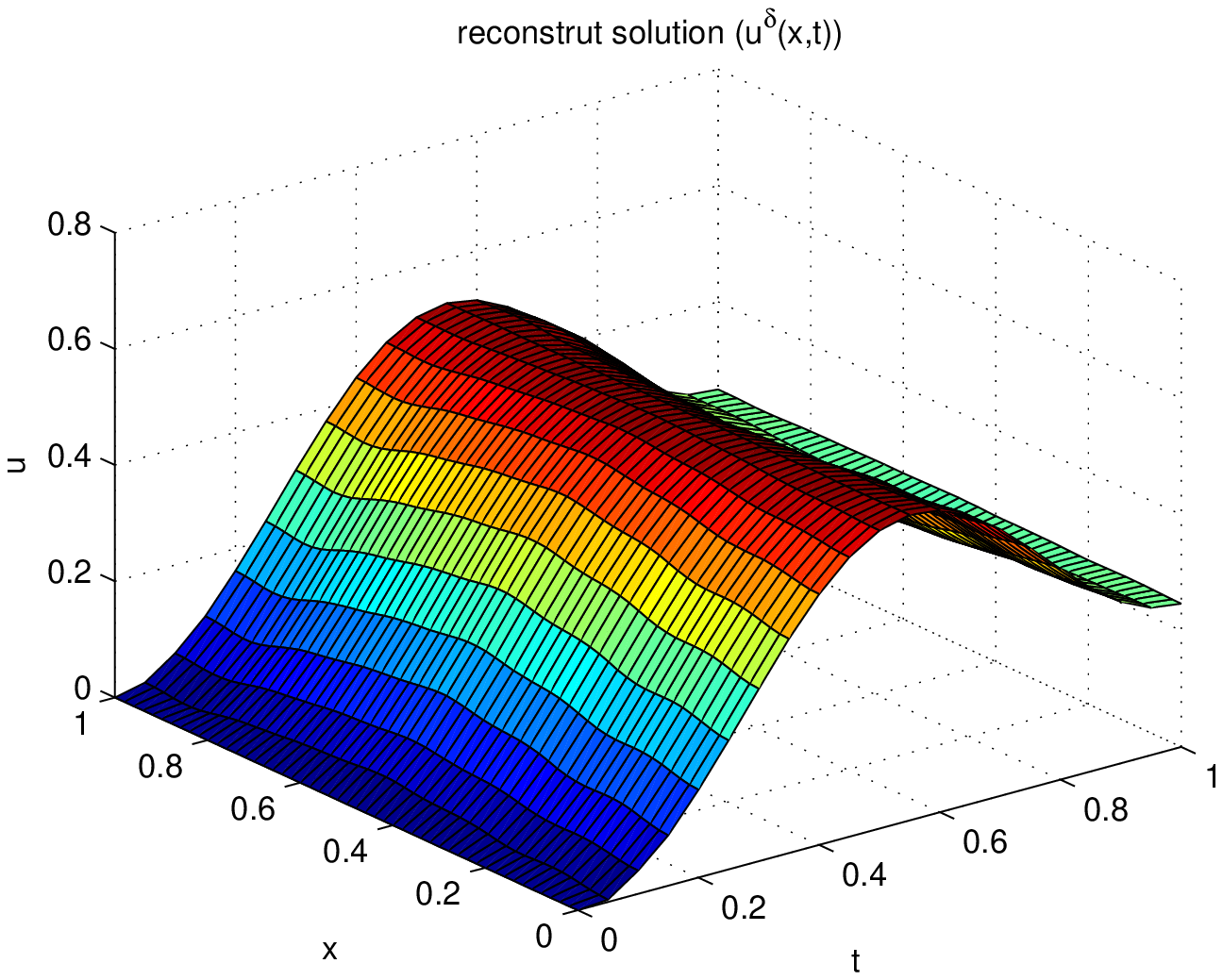}
   \end{center}
    \caption{Comparison of true solution and reconstruct solution. The noise level is $1\%$ and $\vartheta=12$, $p=1/3$, $\sigma=0.2$.
       }\label{fig:5}
\end{figure}

Set $f(x)=\frac{(1-x)x}{|x-\frac{1}{2}|+\frac{1}{2}}$, $\mu_0(x)=0$, $h(t)=5\sin(2\pi t)+1$ and $k=1$. The true solution
and the final data can be calculated by (\ref{eq:302}) (Fig. \ref{fig:4}). Since the final data is nearly a constant function, small
noise level can still produce striking different measurement data comparing with the exact data. Thus here we choose noise level with $1\%$ and $0.1\%$.
Furthermore, we see that $h(t)$ does not satisfy the identically non-positive or non-negative property. However, we
can still get the convergence results since $h(t)$ is not too 'bad'. Fig. \ref{fig:5} shows the reconstructed solution and the exact solution.

\subsection{Example 3}
\begin{figure}[h]
   \begin{center}
         \includegraphics[width=2.8in]{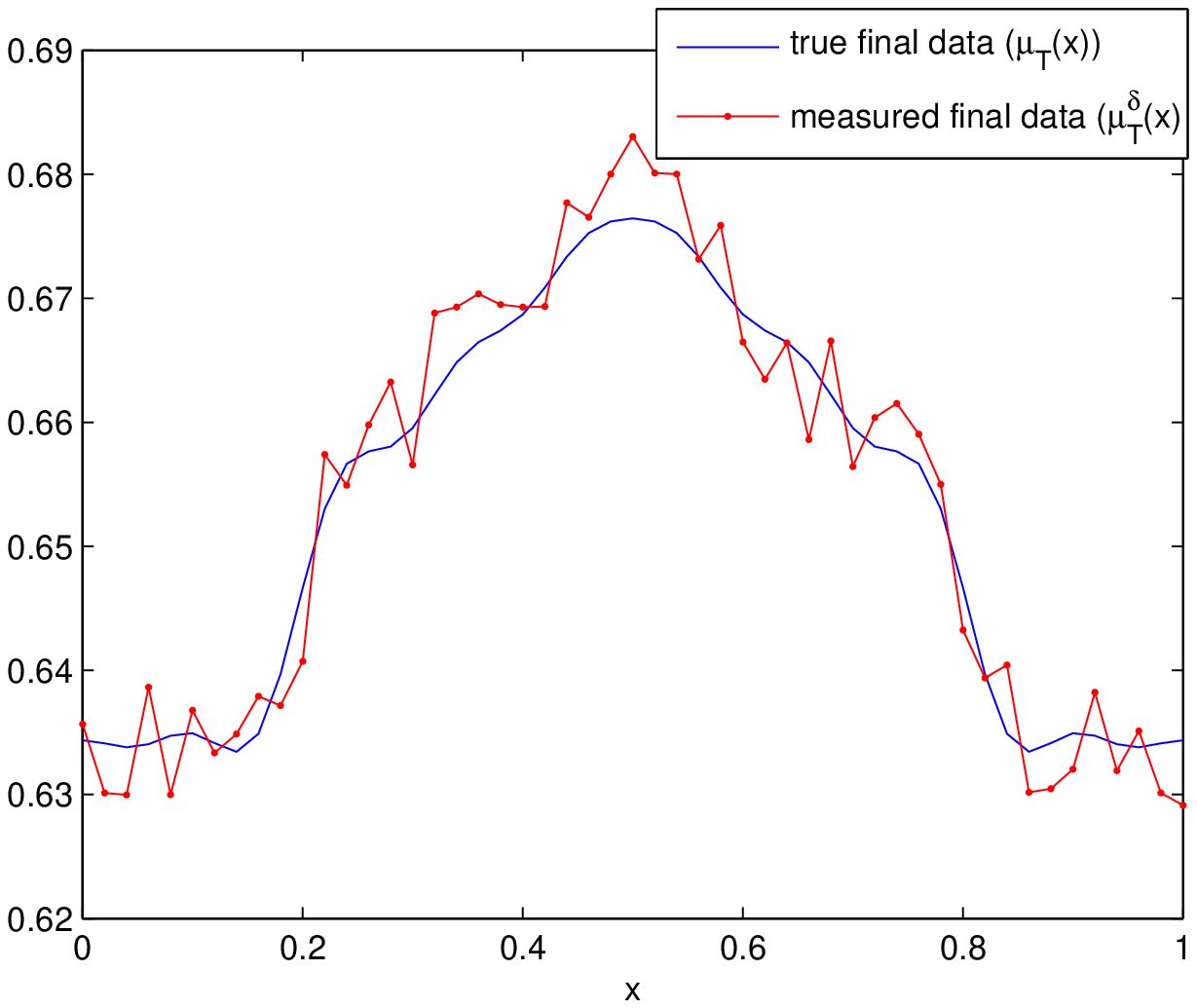}
         \includegraphics[width=2.8in]{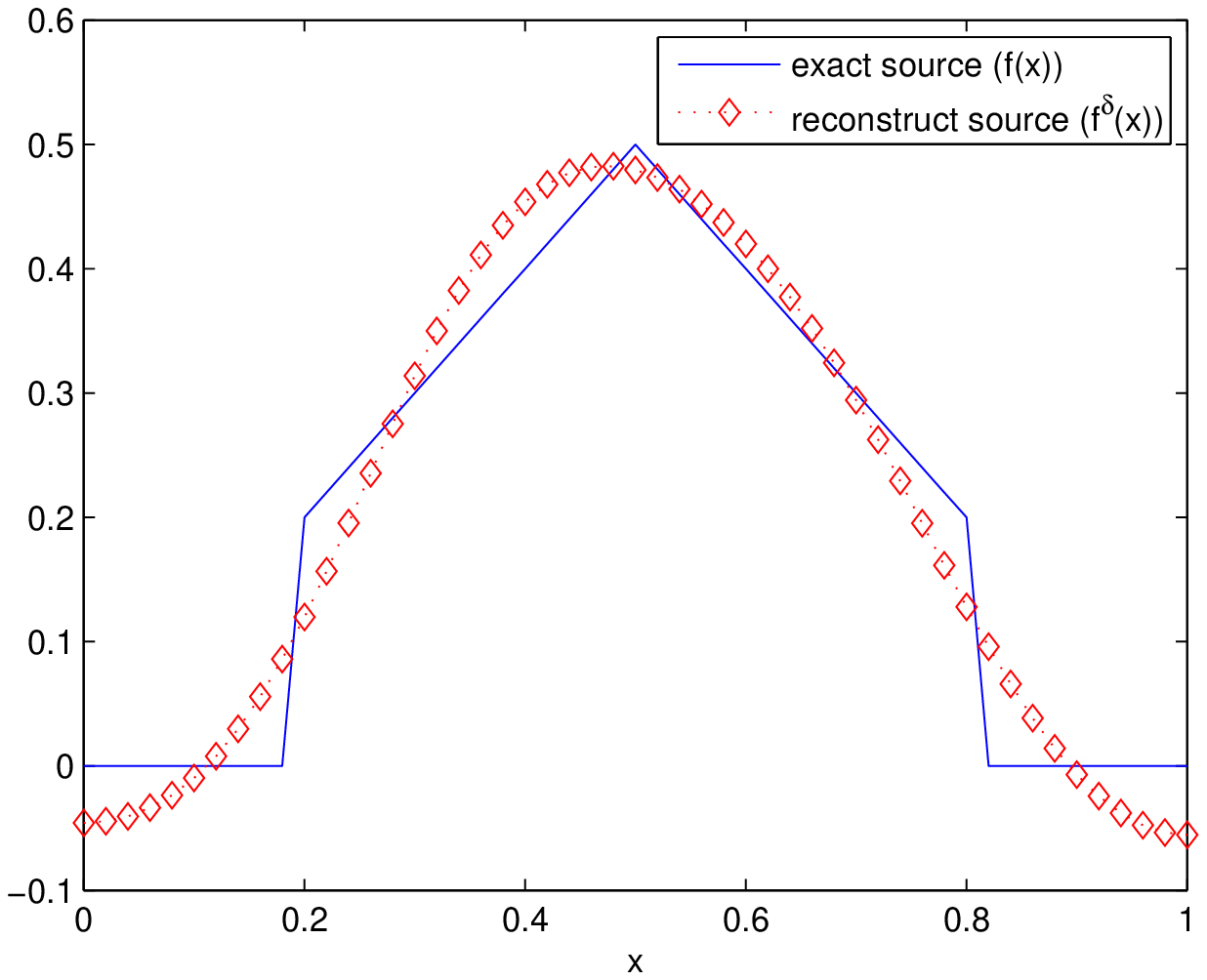}
         \includegraphics[width=2.8in]{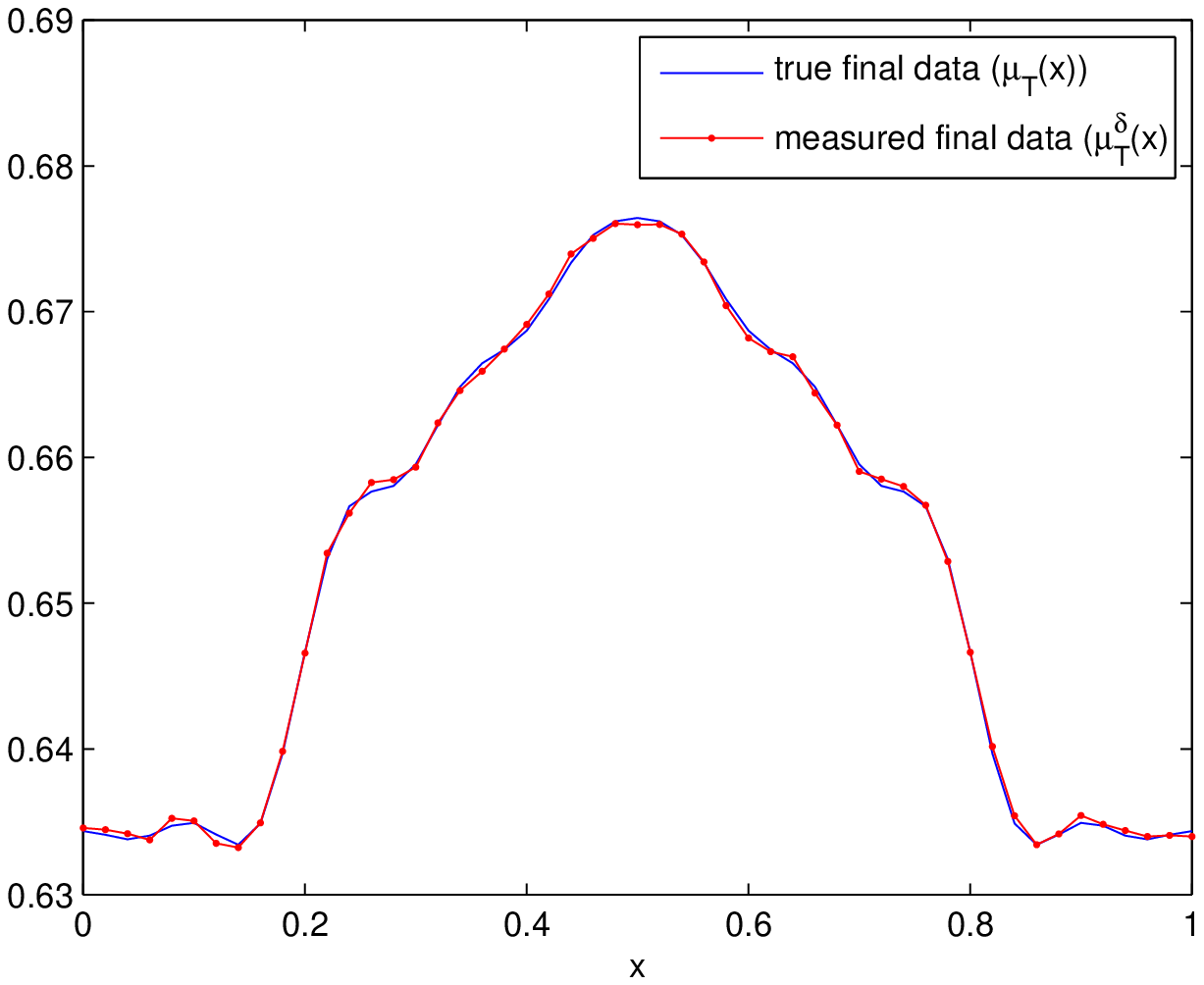}
         \includegraphics[width=2.8in]{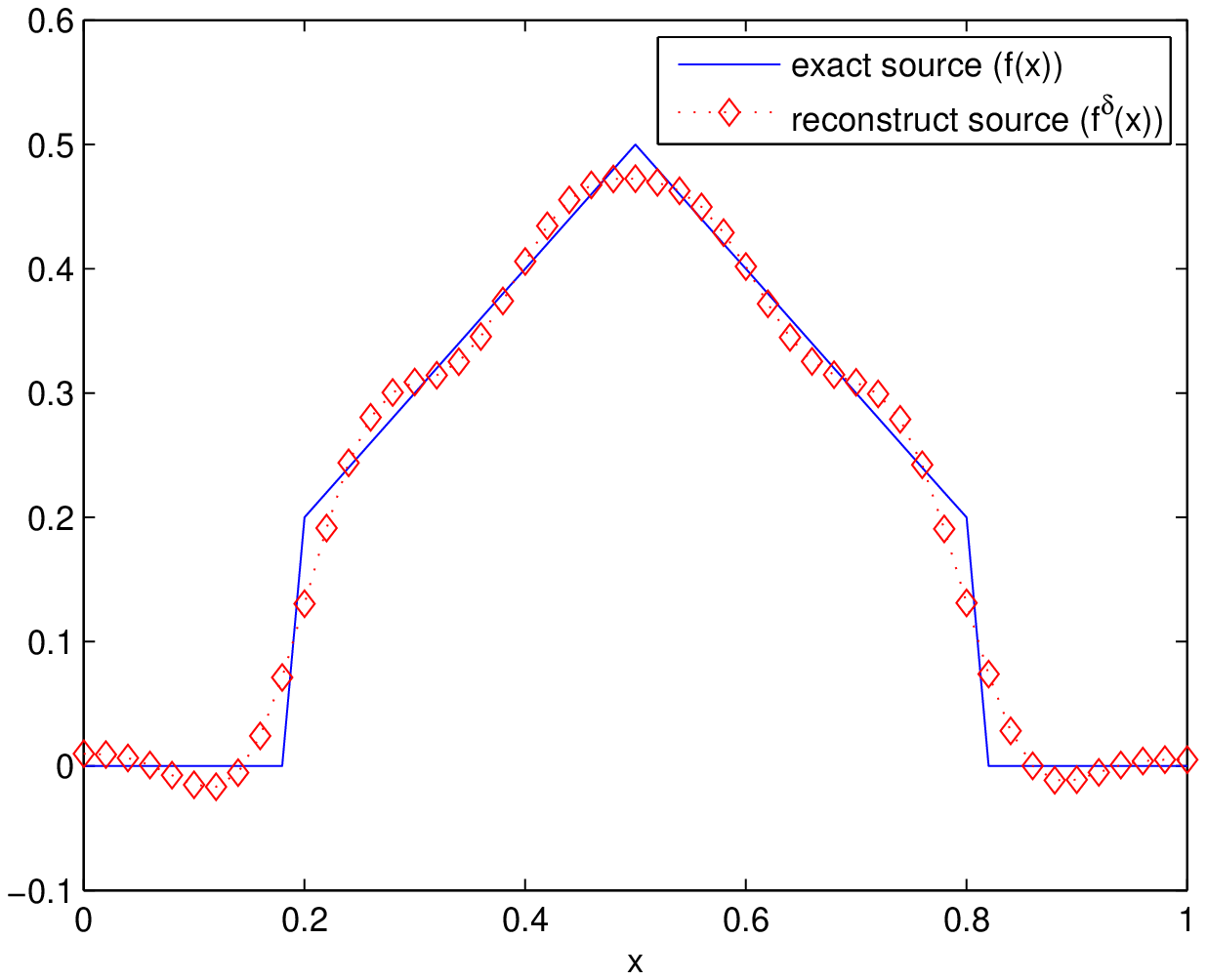}
   \end{center}
    \caption{Comparison of convergence results. The top two with the noise level $1\%$, while the
        bottom two with the noise level $0.1\%$. The left are the final data and corresponding
        measurement error data. The other parameters are $\vartheta=12$, $p=1/3$, $\sigma=0.2$.
       }\label{fig:6}
\end{figure}

\begin{figure}[h]
   \begin{center}
         \includegraphics[width=2.8in]{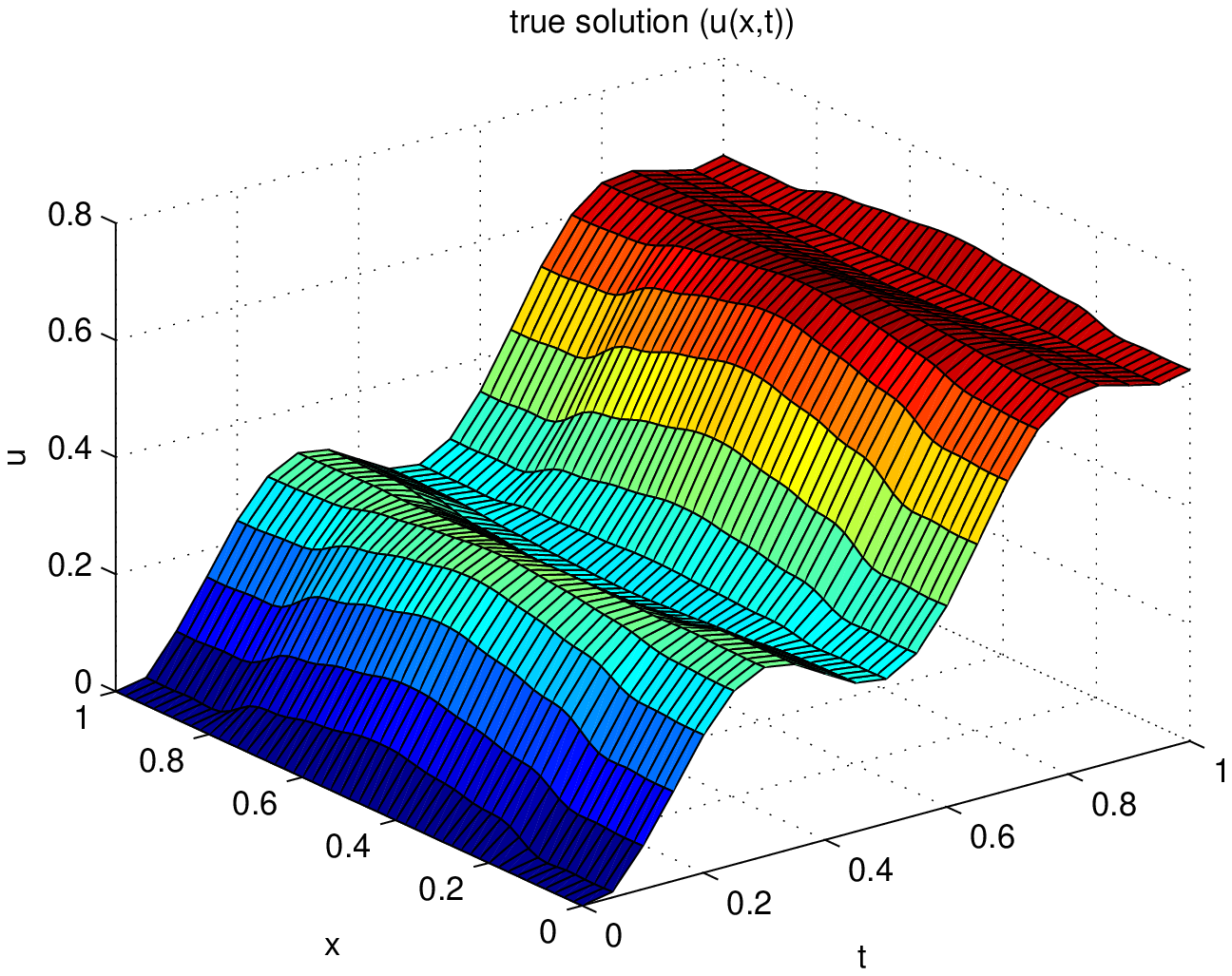}
         \includegraphics[width=2.8in]{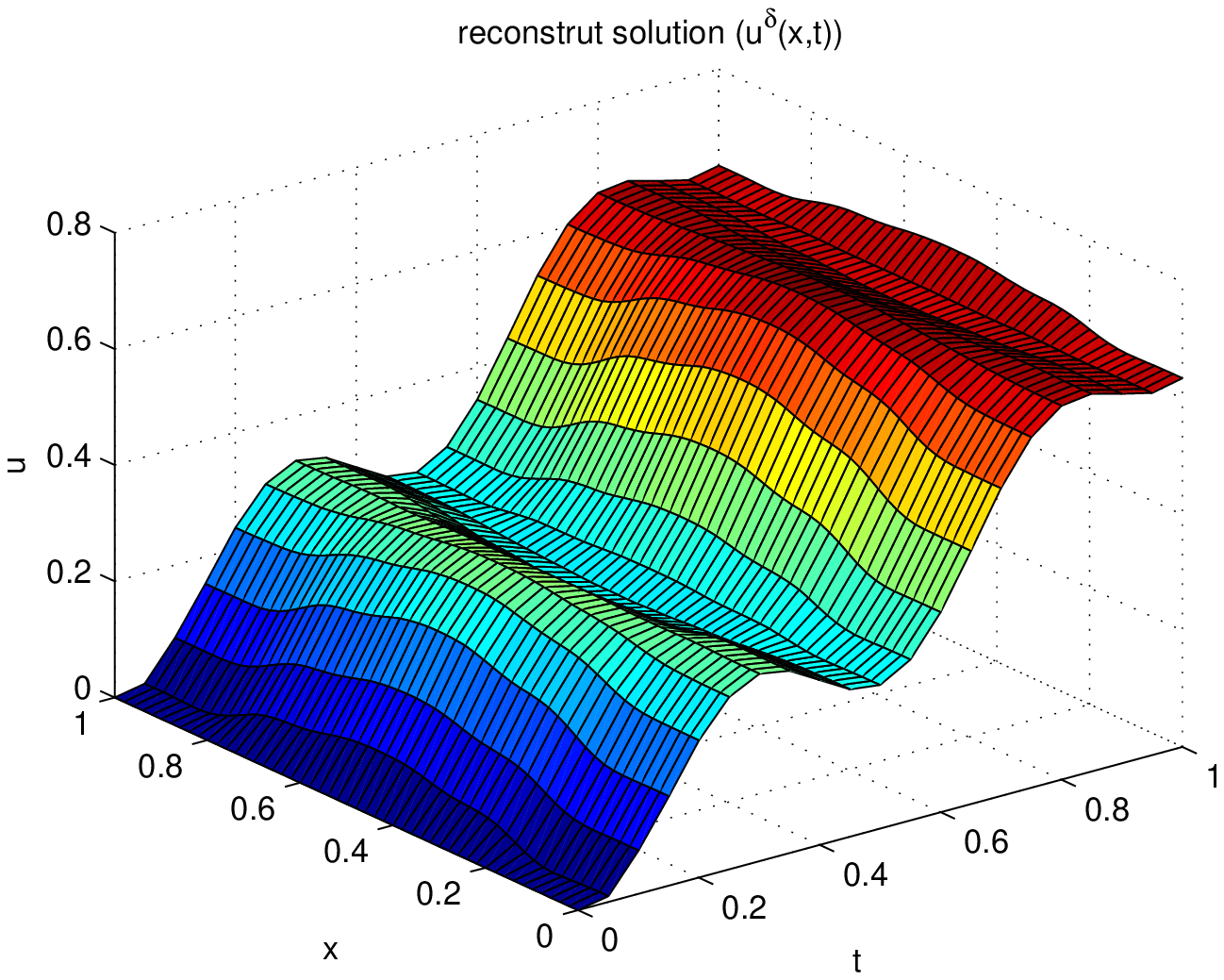}
   \end{center}
    \caption{Comparison of true solution and reconstruct solution. The noise level is $0.1\%$ and $\vartheta=12$, $p=1/3$, $\sigma=0.2$.
       }\label{fig:7}
\end{figure}

Set
$$
f(x)=
\left\{
\begin{array}{cc}
0 & x<0.2 \\
x &  0.2\leq x \leq 0.5 \\
1-x & 0.5 < x \leq 0.8 \\
0 & x > 0.8
\end{array}
\right.
$$
and $\mu_0(x)=0$, $h(t)=e^t+6\sin(4\pi t)+t^2+1$ and $k=1$. The source $f(x)$ has two discontinuous points ($x=0.2, 0.8$) and
one non-differentiable point ($x=0.5$). We see from Fig. \ref{fig:6} that the source can still be reconstructed accurately in a 'smooth' way
when the noise level is not high. Fig. \ref{fig:7} shows the comparison between the true solution and the reconstructed solution.

\section{Conclusions}
In this paper, the numerical methods for reconstruction of source term in both no boundary and Neumann boundary conditions are
presented. The convergence rate has been proved for both {\em a priori} and {\em a posteriori} stopping rules. More importantly,
we show that the solution of the boundary conditions problem has the form of solution for the no boundary problem, which can be
applied for both Neumann and Drichlet boundary conditions. The numerical experiments have shown that the frequency cut-off
technique method applies well for the boundary conditions problem, although for more accurate results we may implement the
iterative methods together with the {\em a posteriori} stopping rule. The numerical methods can be moved parallel to the two
dimensional inverse source problem.

\end{document}